\documentclass[final,leqno,onefignum,onetabnum]{siamltex1213}

\newcommand  \eps{ \varepsilon}
\usepackage{graphicx}
\usepackage{subfig}

\usepackage{amsmath}
\usepackage{textcomp}
\usepackage{enumerate}

\title{  Shilnikov Homoclinic Bifurcation of Mixed-Mode Oscillations }

\author{John Guckenheimer\footnotemark[1]\
\and Ian Lizarraga\footnotemark[2]\ }

\begin{document}
\maketitle
\newcommand{\slugmaster}{%
\slugger{siads}{xxxx}{xx}{x}{x--x}}

\renewcommand{\thefootnote}{\fnsymbol{footnote}}

\footnotetext[1]{Mathematics Department, Cornell University, Ithaca, NY 14853}
\footnotetext[2]{Center for Applied Mathematics, Cornell University, Ithaca, NY 
14853}

\renewcommand{\thefootnote}{\arabic{footnote}}

\begin{abstract}
The Koper model is a three-dimensional vector field that was developed to study 
complex electrochemical oscillations arising in a diffusion process.
Koper and Gaspard described paradoxical dynamics in the model: they discovered 
complicated, chaotic behavior consistent with a homoclinic orbit of Shil'nikov 
type, but were unable to locate the orbit itself. The Koper model has since 
served as a prototype to study the emergence of mixed-mode oscillations (MMOs) 
in slow-fast systems, but only in this paper is the existence of these elusive 
homoclinic orbits established. They are found first in a larger family that has 
been used to study singular Hopf bifurcation in multiple time scale systems 
with two slow variables and one fast variable. A curve of parameters with 
homoclinic orbits in this larger family is obtained by continuation and shown 
to cross the submanifold of the Koper system. The strategy used to compute the 
homoclinic orbits is based upon systematic investigation of intersections of 
invariant manifolds in this system with multiple time scales. Both canards and 
folded singularities are multiple time scale phenomena encountered in the 
analysis. 
Suitably chosen cross-sections and return maps illustrate the complexity of the 
resulting MMOs and yield a modified geometric model from the one Shil'nikov 
used to study spiraling homoclinic bifurcations.
\end{abstract}


\begin{AMS}\end{AMS}

\pagestyle{myheadings}
\thispagestyle{plain}
\markboth{JOHN GUCKENHEIMER AND IAN LIZARRAGA}{SHILNIKOV HOMOCLINICS IN THE 
KOPER MODEL}

\section{Introduction} \label{sec:introduction}

In 1992, Marc Koper and Pierre Gaspard introduced a three-dimensional model to 
analyze an electrochemical diffusion problem, in which layer concentrations of 
electrolytic solutions fluctuate nonlinearly at an electrode 
\cite{koper1991jphysrev, koper1992jchemphys}. They sought to model mixed-mode 
oscillations (hereafter MMOs) arising in a wide variety of electrochemical 
systems. Their analysis revealed a host of complicated dynamics, including 
windows of period-doubling bifurcations, Hopf bifurcations, and complex Farey 
sequences of MMO signatures.

They also found regions in the parameter space where the equilibrium point of 
the system satisfies the Shilnikov condition.\footnote[1]{Let the linearization 
of a three-dimensional vector field at an equilibrium point $p$ have 
eigenvalues 
$\rho \pm i \omega$ and $\lambda$, where $\rho, ~\omega,$ and $\lambda$ are all 
real. Then $p$ satisfies the {\it Shilnikov condition} if 
$\rho \lambda < 0$ and $|\rho/\lambda| < 1$.} Within these regions, they 
observed that trajectories  repeatedly come close to the fixed point, and 
return 
maps strongly suggest chaotic motion consistent with a Shilnikov homoclinic 
bifurcation.  However, they were unable to locate a genuine homoclinic orbit to 
account for this behavior, so it was catalogued as a {\it near homoclinic 
scenario}. In such a scenario, complex and chaotic MMOs could suddenly 
arise---as if from a homoclinic bifurcation---but without the existence of the 
homoclinic orbit to serve as an organizing center.

Nevertheless, the Koper model has emerged as a paradigm in studies of slow-fast 
systems containing MMOs~\cite{braaksma1998jnonsci,desroches2012siamreview}. As 
its four parameters are varied, local mechanisms such as folded nodes 
\cite{guckenheimer2008chaos, wechselberger2005siads} and singular Hopf 
bifurcation \cite{baer1986sjaa, braaksma1998jnonsci, guckenheimer2008siads, 
guckenheimer2012siads, guckenheimer2012dcdsa} generate small-amplitude 
oscillations. A global return mechanism allows for repeated reinjection into 
the regions containing these local objects. The interplay of these local and 
global mechanisms gives rise to sequences of large and small oscillations 
characterized by signatures that count the numbers of consecutive small and 
large oscillations. Shilnikov homoclinic orbits are limits of families of MMOs 
with an unbounded number of small oscillations in their signatures. When the 
Shilnikov condition is satisfied, they also guarantee the existence of chaotic 
invariant sets whose presence in the Belousov-Zhabotinsky reaction was 
controversial for several years.

The homoclinic orbits that could explain Koper's original observations have 
remained elusive. This paper describes their first successful detection. 
Multiple  timescales make numerical study of these homoclinic orbits quite 
delicate; on the other hand, the presence of {\it slow manifolds} allows us to 
analyze many aspects of the system with low-dimensional maps. We find the 
homoclinic orbits by exhibiting multiple time scale phenomena such as canards 
and folded singularities. We first locate such orbits in a five-parameter family 
of 
vector fields used to explore the dynamics of singular Hopf bifurcation. After 
an affine coordinate change, the Koper model is a four-parameter subfamily. 
Shooting methods that compute trajectories between carefully chosen 
cross-sections cope with the numerical instability resulting from the singular 
behavior of the equations. Following identification of a homoclinic orbit in 
the 
larger family, a continuation algorithm is used to track a curve on the 
codimension-one 
manifold of spiraling homoclinic orbits in parameter space. This manifold 
intersects the parametric submanifold 
corresponding to the Koper model, locating the homoclinic orbit that is the 
target of our search.

This paper is a numerical investigation of \emph{the} Koper model. The results 
are not rigorous, so a formal ``theorem-proof'' style is inappropriate for this 
discussion. Following the numerical part of the paper, we abstract our 
reasoning to present a geometric model for the slow-fast decomposition of the 
homoclinic orbits of the Koper model. This geometric model produces a list of 
properties that we use to prove the existence of a homoclinic orbit in a 
slow-fast system. We think these properties may be amenable to verification in 
the Koper model through the use of interval arithmetic.

\section{The Koper model as a slow-fast system} \label{sec:koper}

We consider vector fields of the form
\begin{eqnarray} 
\varepsilon \dot{x} &=& f(x,y), \nonumber\\
\dot{y} &=& g(x,y),
\end{eqnarray} \label{eqn1}
where $x \in R^m$, $y \in R^n$, and the functions $f$ and $g$ are smooth. In 
this paper, $f$ and $g$ are polynomials. {\it Slow-fast}  vector fields are 
those where $\varepsilon \ll 1$. In this case, $x$ is the fast variable and $y$ 
is the slow-variable. 

The Koper model is a frequently studied example 
and the subject of this paper. It is defined by 
\begin{eqnarray}
\varepsilon_1 \dot{u} &=& kv - u^3 + 3u - \lambda, \nonumber\\
\dot{v} &=&  u - 2v + w, \label{eq:koper}\\
\dot{w} &=& \varepsilon_2 (v - w), \nonumber
\end{eqnarray}
and has two slow variables $v$ and $w$ and one fast variable $u$ when the 
parameters $\varepsilon_1 \ll 1$ and $\varepsilon_2 = 1$. The terms $k$ and 
$\lambda$ are additional parameters. We assume throughout this paper that 
$\varepsilon_2 = 1$ without further comment.

The study of slow-fast vector fields has advanced rapidly in recent years with 
specific focus on systems having two slow variables and one fast variable. We 
recall relevant aspects of the theory that bear on the Koper model. The review 
paper \cite{desroches2012siamreview} includes an extensive discussion that 
provides additional information, especially about mixed mode oscillations. 

The set of points defined by $ C = \{f = 0\}$ in system \eqref{eqn1} is called 
the {\it critical manifold}. 
Fenichel \cite{fenichel1972indunivmathj} proved the existence of locally 
invariant slow manifolds near regions of  $C$ where $D_x f$ is hyperbolic. 
Trajectories on the slow manifolds are approximated by trajectories of the {\it 
reduced system} on $C$ defined by
\begin{eqnarray}
\dot{y} &=& g(h(y),y,0),
\end{eqnarray}
where $h$ is defined implicitly by $f(h(y),y,0) = 0$.  While the slow manifolds 
are not unique, compact portions are exponentially close to each other: their 
distances from each other are $O(\exp(-c/\varepsilon))$ as $\eps \to 0$. We 
often refer to `the' slow manifold in statements where the choice of slow 
manifold does not matter. The points $x \in C$ where $D_x f$ is singular are 
called {\it fold points}. For the Koper model, $C$ is the zero set of $kv - u^3 
+ 3u - \lambda$ and the fold curve consists of the points on $C$ with 
$\frac{3}{k}(u^2-1) = 0$. The reduced system is given by
\begin{eqnarray}
\frac{3}{k}(u^2-1) \dot{u} &=&  u - 2v + w, \label{eq:koperslow}\\
\dot{w} &=& (v - w), \nonumber
\end{eqnarray}
where $v= \frac{1}{k}(u^3 - 3u + \lambda)$.

The reduced system is a \emph{differential algebraic equation} that is the main 
component of the \emph{singular limit} of the ``full'' system \eqref{eq:koper} 
as $\eps_1 \to 0$. The other component of the singular limit comes from the 
fast part of the full system. The fold curve divides the critical manifold into 
attracting and repelling sheets where $u^2 > 1$ and $u^2 < 1$ respectively. 
When trajectories of the full system reach the vicinity of the fold curves, 
they turn and then flow close to the fast direction to a small neighborhood of 
the 
attracting slow manifold. In the singular limit, one has trajectories on the 
attracting and repelling sheets of the critical manifold that meet at the fold 
curve separating the sheets. The limiting behavior of the full system from 
these impasse points is to \emph{jump} along the fast direction stopping at the 
first intersection with an attracting sheet of $C$. The addition of fast jumps 
makes the singular limit a \emph{hybrid} system whose trajectories (called 
\emph{candidates}) consist of concatenations of slow segments that solve the 
reduced system and fast jumps parallel to the fast direction. As we describe 
below, fast jumps may occur at points of $C$ which are not fold points.

Rescaling time of the reduced system produces a 
\emph{desingularized} reduced system that extends to the fold curve. In the 
case 
of the Koper model, the desingularized reduced system is
\begin{eqnarray}
\dot{u} &=&  u - 2v + w, \label{eq:koperslowdesing}\\
\dot{w} &=& \frac{3}{k}(u^2-1) (v - w), \nonumber
\end{eqnarray}
where $v= \frac{1}{k}(u^3 - 3u + \lambda)$.
Discrete time jumps of trajectories that reach the fold curve of 
\eqref{eq:koperslowdesing} are added to the reduced model to reflect the 
limiting behavior of the full system. However, there are exceptional  
trajectories of the full system that flow from an attracting slow manifold to a 
repelling slow manifold without making a jump at the fold curve. The singular 
limits of these trajectories flow to equilibrium points of 
\eqref{eq:koperslowdesing} which lie on its fold curves. Since the time 
rescaling of desingularization reverses the direction of time on the 
repelling sheet of $C$ where $u^2 < 1$, we follow trajectories of 
\eqref{eq:koperslowdesing} backwards in time ``through'' the singularity. At 
any location along such a backward trajectory, a jump to one of the attracting 
sheets of the slow manifold is allowed along the fast direction. This makes the 
reduced system multi-valued and hence still more complicated. See~\cite{GHW03b} 
for 
further discussion of this construction of a hybrid reduced model in the context 
of the forced van der Pol system. 

\begin{figure}[!ht]
    \centering
    \subfloat[ ]{{\includegraphics[width=3in, height=1.7in]{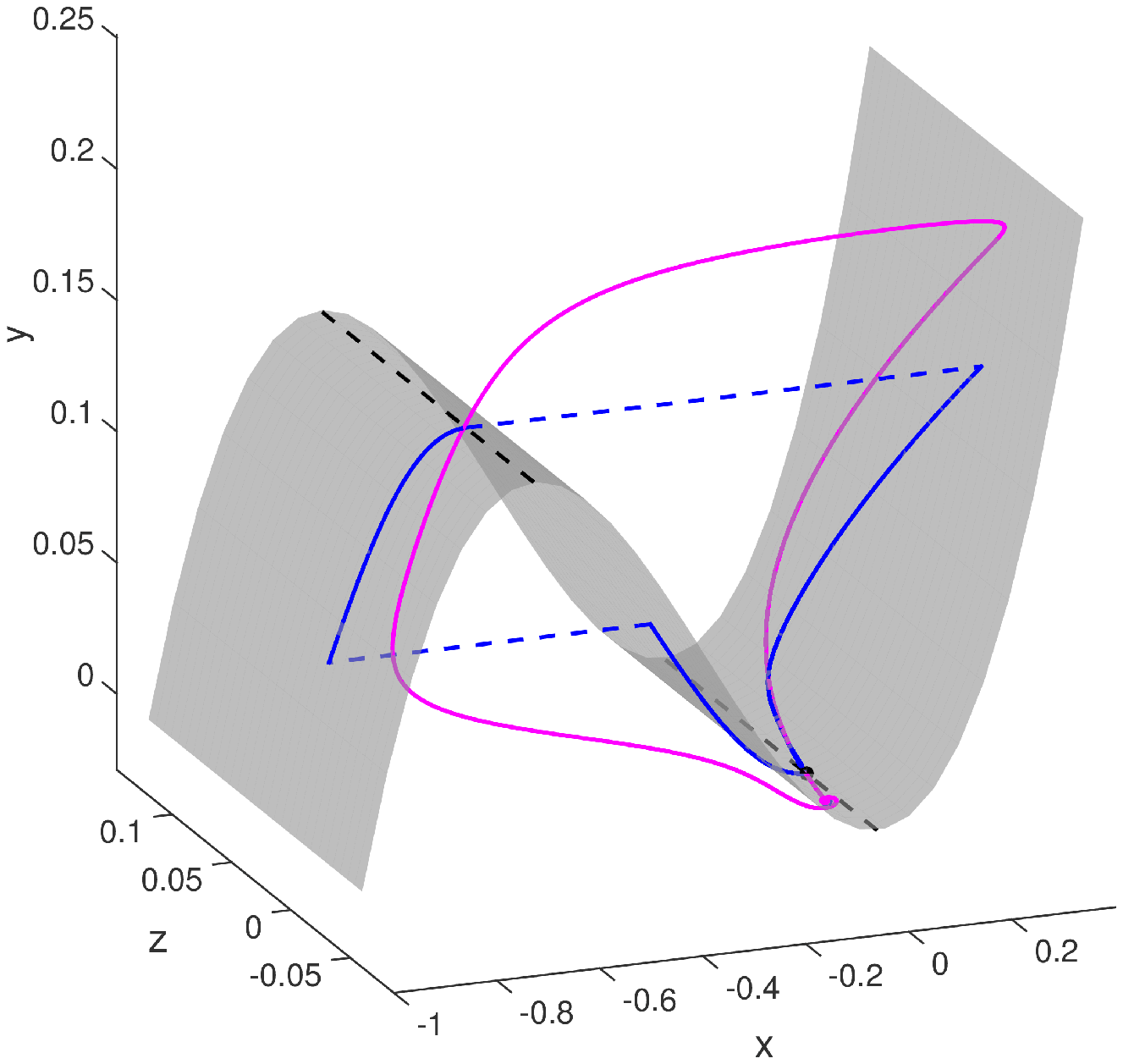} }} 
    \subfloat[ ]{{\includegraphics[width=3in, height=1.7in]{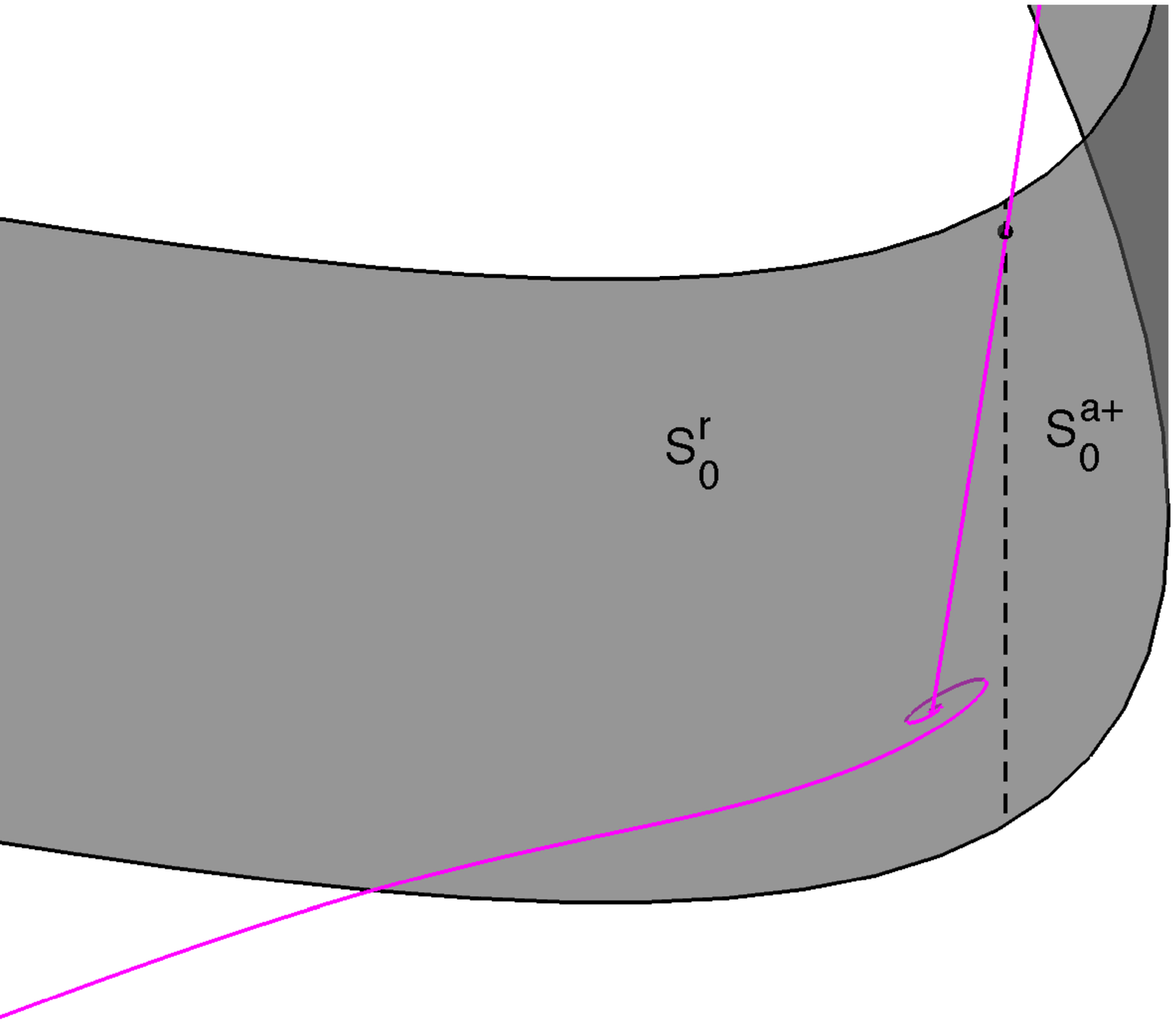} }}
    \caption{(a) Shilnikov homoclinic orbit (magenta curve) found in the 
Koper model, defined by \eqref{eq:shnf}, alongside a typical singular 
homoclinic 
orbit defined by trajectories (blue solid curves) lying on the critical 
manifold 
(gray surface) of the system, concatenated with instantaneous jumps (blue 
dashed curves) to different branches of the critical manifold. Fold lines are 
given by black dashed curves. 
Section \ref{sec:singularorbits} discusses families of singular homoclinic 
orbits in detail.
(b) A cartoon of the local 
dynamics. The homoclinic orbit departs along the two-dimensional unstable 
manifold of the equilibrium and returns by passing 
very close to a folded singularity of the singular system,
given by the black dot on the fold line separating the sheets $S^r_0$ and 
$S^{a+}_0$ of the critical manifold. These sheets perturb to attracting and 
repelling slow manifolds for small values of the singular perturbation 
parameter $\eps$.}
    \label{fig:descriptive}
\end{figure}

The desingularized reduced system of a fast-slow system can have two types of 
equilibrium points. Equilibria of the ``full'' system are retained as 
equilibria of the desingularized reduced system. Additionally, points on a fold 
curve may become equilibria of the 
desingularized reduced system because the time rescaling factor vanishes at 
these points. For the Koper model, these \emph{folded 
singularities} are solutions of the equations $kv - u^3 + 3u - \lambda$, $u - 
2v 
+ w = 0$ and $u^2 - 1 = 0$. (The third of these equations replaces the 
equation $\dot{z} = 0$ for an actual equilibrium.) The folded singularities 
are further classified as {\it folded nodes, saddles, saddle-nodes}, and {\it 
foci} depending 
on their type as equilibria of the desingularized reduced system. The folded 
singularities mark transitions along the fold curve where the trajectories of 
the reduced system (without desingularization) approach the fold curve or flow 
away from it. They are also close to places where the full system might have 
trajectories that contain segments on an attracting slow manifold that proceed 
to follow the repelling slow manifold without an intervening jump. Trajectory 
segments that flow for time $O(1)$ in the slow time scale along a repelling 
slow manifold  of a slow-fast systems are called \emph{canards}. They have 
been the subject of intense study since they were found in periodic orbits of 
the forced van der Pol equation \cite{BenoitCallotDienerDiener,diener}. They 
also play an important role in the dynamics of the Koper model and are part of 
the homoclinic orbits we find. See Figure \ref{fig:descriptive}.

The dynamics of the full system close to folded singularities is significantly 
more complicated than suggested by the reduced system. Benoit~\cite{benoit} 
observed that trajectories of the full system near folded nodes of the reduced 
system can make small amplitude oscillations in directions that involve both 
fast and slow variables. The oscillations are manifest in twisting of the 
attracting and repelling slow manifolds, as visualized by Desroches 
\cite{desroches2012siamreview}. The maximum number of small oscillations of 
trajectories passing through this region is bounded, with a bound determined 
from a normal form for the folded node \cite{benoit,wechselberger2005siads}. 
In the Koper model and other slow-fast systems, small oscillations associated 
with folded nodes together with ``global returns'' produce \emph{mixed mode 
oscillations} \cite{desroches2012siamreview}.

Our goal in this paper is to locate a homoclinic orbit of the Koper model when 
its equilibrium is a saddle-focus. The initial portion of such an orbit 
consists of small amplitude oscillations as the homoclinic orbit spirals 
away from the equilibrium along its unstable manifold (see Fig. 
\ref{fig:descriptive}). In single time scale 
systems, both shooting and boundary value methods can be used to find the 
homoclinic orbit. With both types of algorithms, as well as in theoretical 
analysis of the dynamics near the homoclinic orbit, trajectories are decomposed 
into a local portion in a neighborhood of the equilibrium and a 
global portion that lies outside that neighborhood. Linearization or higher 
order approximations are used to identify the local stable and unstable 
manifolds, and then numerical methods are used to find the global 
return. 

Because the homoclinic bifurcation is a codimension one phenomenon, an 
``active'' parameter is used in locating the bifurcation. Extending the phase 
space  to include the active parameter, the homoclinic orbit becomes a 
transverse intersection of the families of stable and unstable manifolds of the 
curve of equilibria in the extended system. Shooting methods locate the global 
return by computing a two parameter family of initial value problems whose 
initial conditions depend upon the active parameter and a coordinate that 
parameterizes trajectories in the unstable manifold. The intersections of these 
trajectories with a cross-section of the local stable manifolds forms a surface 
while the family of local stable manifolds intersects the cross-section in a 
curve. The cross-section is three dimensional in the extended phase space, so 
the intersections of the stable and unstable manifolds can cross transversally. 

Straightforward implementation of this strategy fails in the Koper model, as 
Koper discovered in his original investigation. The slow-fast 
structure of the problem makes portions of the global return trajectories so 
sensitive that continuous dependence upon initial conditions appears 
to fail for numerical solutions of initial value problems. This apparent 
discontinuity of trajectories has been documented in the ``canard 
explosions'' of the van der Pol equation with constant forcing~\cite{GHW}. The 
crux of our analysis of the Koper homoclinic model is to decompose the global 
returns into segments, each of which can be reliably computed with either 
forward or backward integration. This strategy has worked with other 
examples
like the canard explosion \cite{GL2007} of the van der Pol system with 
constant forcing, and it proves successful again here.

When the Koper model equilibrium is a saddle-focus, its two dimensional
unstable manifold is a source of small amplitude oscillations. 
The number of small amplitude oscillations is unbounded as 
trajectories leave the equilibrium. However, their magnitude changes 
quickly 
unless the ratio of the real and complex parts of the complex eigenvalues at 
the equilibrium is small. Parameter regions close to a Hopf bifurcation 
yield small ratios, so choosing parameters in such a region is part of our 
strategy for finding homoclinic orbits. Hopf bifurcation in slow-fast systems 
takes several forms 
that depend upon whether the two dimensional subspace of center directions lie 
in fast directions, slow directions, or a mixture of the two. The third case 
is called \emph{singular Hopf bifurcation}, and it is the case that occurs in 
the Koper model. In systems with two slow variables, singular Hopf bifurcation 
is closely associated with a phenomenon in the desingularized reduced system 
called \emph{folded saddle-node, type II} (FSNII)\cite{wechselberger2005siads}. 
This occurs in the desingularized reduced system when an equilibrium (that is 
not a folded singularity) crosses a fold curve as a parameter is varied. As the 
equilibrium crosses the fold curve, a folded singularity passes through the 
same point, and its type switches between a fold node and a folded saddle. 
Singular Hopf bifurcations are found at a distance $O(\eps)$ from the FSNII 
point.

The dynamics of a model system with an FSNII bifurcation have been analyzed by 
Guckenheimer and Meerkamp~\cite{guckenheimer2012siads}. The model they 
studied is obtained 
by truncating a Taylor expansion at the FSNII bifurcation: it can be viewed 
informally as a normal form for this problem and for singular Hopf bifurcation. 
Moreover, this normal form is closely related to the Koper model. With 
the addition of a single cubic term and an appropriate affine coordinate 
change, 
it contains the Koper system as a subfamily. We make use of this relationship 
below, first finding a homoclinic orbit of the larger family and then 
continuing 
it to parameters that lie in the (transformed) Koper family. The FSNII normal 
form has many different types of bifurcations and small amplitude chaotic 
behavior is possible~\cite{guckenheimer2012siads}. Complicated small scale 
dynamics of this system raise the possibility of homoclinic orbits that are 
more 
complicated than the ones we exhibit in this paper. 

A key aspect of the FSNII dynamics relevant to our search for homoclinic 
orbits in the Koper model is the intersection of the unstable manifold of the 
equilibrium with the repelling slow manifold. Since both of these manifolds are 
two dimensional, we expect them to intersect along isolated trajectories. The 
homoclinic orbits we find contain segments close to such intersections. As 
trajectories in 
the unstable manifold emerge from the equilibrium point, they may follow the 
repelling slow manifold for some distance, producing a canard that can then 
jump 
away from the repelling manifold. Part of our computational task is to identify 
the jump point that yields the homoclinic orbit.

By definition, the homoclinic orbit contains a branch of the stable manifold 
of the equilibrium. We find that this branch lies close to the attracting slow 
manifold for a substantial distance. A second key part of our 
computations stems from 
the observation that the stable manifold of the equilibrium crosses the 
attracting slow manifold as a parameter is varied and that the homoclinic 
bifurcation lies exponentially close to this parameter value. Finding this 
crossing is delicate because the stable manifold passes through a region 
where extensions of the normally hyperbolic attracting 
manifold twist. Most of the 
studies of mixed mode oscillations of the Koper model to date are based upon 
local analysis of oscillations close to perturbations of a folded node or near 
the equilibrium. Here, we need to study trajectories that interact with both 
the equilibrium and a twisting, attracting slow manifold, a situation that 
remains poorly 
understood. Our 
homoclinic orbit lies close to such trajectories, and further analysis of their 
dynamics is likely to produce interesting results.

We close this section with three remarks:
\begin{enumerate}
 \item 
Folded singularities are defined for the reduced system that is the 
singular 
limit of a slow-fast system. In the full system with time scale parameter 
$\eps>0$, these entities are no longer defined. However, it is useful to 
identify points that are located where their influence is manifest. We use the 
term \emph{twist region} to describe sets in the phase space of the full 
system where an (extended) attracting slow manifold twists around a canard. 
Numerically, we locate approximations to these points as folded nodes of the 
reduced system obtained by setting $\eps = 0$.
\item
We consider only parameter values of the Koper model for which its 
equilibrium point $p_{eq}$ has a pair of complex conjugate eigenvalues $\rho 
\pm i \omega $ and one real eigenvalue $\lambda$ with $\lambda < 0 < \rho$.
This constraint forces the equilibrium to be within a distance $O(\eps)$
of the folded singularity of the singular limit. This requirement is explained 
further in Section~\ref{sec:shnf} below. When we study the singular limit 
of the system, we do so along parameter curves along which the equilbrium 
remains a saddle-focus.
\item
Numerical studies of slow-fast systems always rely upon particular choices of 
the time scale parameter $\eps > 0$, while the theory focuses upon dynamics 
that 
is present for ``small enough $\eps$.'' It is rare that one can verify that a 
particular choice of $\eps$ is indeed small enough to fall within the scope of 
particular theorems. Nonetheless, the theory provides a guide to expected 
behavior in the numerics, explaining observations that would otherwise appear 
anomalous. These studies are based upon the presumption that $\eps$ is 
sufficiently small and test this presumption by examining the dynamics for 
different values of $\eps$. 
However, like the theory, the problem of connecting 
the behavior observed at specific values of $\eps > 0$ with the singular limit 
is nontrivial. Our approach here is to show that we find behavior 
\emph{consistent} with that described by perturbations from the singular limit.
\end{enumerate}

\section{The singular Hopf extension of the Koper model} \label{sec:shnf}

We now introduce a family of vector fields with two slow variables and one fast 
variable that contains a singular Hopf bifurcation \cite{baer1986sjaa, 
braaksma1998jnonsci, guckenheimer2008siads, guckenheimer2012siads, 
guckenheimer2012dcdsa}. It is given by

\begin{eqnarray}
\varepsilon \dot{x} &=& y - x^3 - x^2, \nonumber\\
\dot{y} &=& z - x, \label{eq:shnf}\\
\dot{z} &=& - \nu - a x - b y - c z, \nonumber
\end{eqnarray}
\noindent
where $\varepsilon, \nu, a, b,$ and $c$ are parameters, $x$  is the fast 
variable, and $y$ and $z$ are the slow variables. We denote $\alpha = 
(\varepsilon, \nu, a, b, c)$ and define $P$ to be the five dimensional 
space of parameters 
$\alpha$. The critical manifold is the $S$-shaped cubic surface $\{y = 
x^3 + x^2\}$ with two {\it fold lines} at $L_{0}$ defined as $\{x = 0\}$ and 
$L_{-2/3}$ defined as $\{x = - 2/3\}$. In our numerical investigations, we 
set $\varepsilon = 0.01$ unless otherwise noted.

``The'' slow manifold $S$ has sheets $S^{a-}_{\varepsilon}$, 
$S^r_{\varepsilon}$ and  $S^{a+}_{\varepsilon}$ that lie close to the sheets 
of the critical manifold
$C$ defined by $C \cap \{ x < -2/3\}$, $C \cap \{-2/3 < x < 0\}$ and $C \cap \{ 
0 < x \}$. Away from the fold lines, forward trajectories are attracted to 
$S^{a\pm}_{\varepsilon}$ and repelled from $S^r_{\varepsilon}$ at fast 
exponential rates (see for eg. \cite{jones1995lecnotesmath} for a derivation of 
estimates using the Fenichel normal form). 

After the affine coordinate change defined by 
$(x,y,z) = ((u-1)/3,(k v - \lambda 
+2)/27,2v-w-1)/3)$~\cite{desroches2012siamreview}, scaling time by $-k/9$ and 
the substitutions 
$$(\eps,a,b,c,\nu) = (-k \eps_1/81, 18/k, 81 \eps_2/k^2,-9(\eps_2+2)/k,
(3\lambda - 6  -3 k)\eps_2/k^2),$$ 
the Koper model becomes a parametric subfamily of \eqref{eq:shnf}, with 
parameters satisfying the equation
\begin{eqnarray}
2b + a(a + c) = 0. \label{kopeqn}
\end{eqnarray}
Note that the above parametric equation corrects a sign error in Desroches et 
al.\cite{desroches2012siamreview}. We work henceforth with the Koper model in 
the form given by \eqref{eq:shnf}.

The corresponding desingularized slow flow of the system \eqref{eq:shnf} is 
given by
\begin{eqnarray}
\dot{x} &=& z - x\nonumber\\
\dot{z} &=& -(2x + 3x^2) (\nu + ax + b(x^2 + x^3) + cx).\label{eq:shnfslow}
\end{eqnarray}
In analogy to \eqref{eq:koperslowdesing}, this two-dimensional system is 
derived 
by rescaling time with the term $(2x + 3x^2)$. The origin is always an
equilibrium of \eqref{eq:shnfslow}, so $n = 
(0,0,0) \in L_0$ is a folded singularity of \eqref{eq:shnf}.

We consider only parameter sets of \eqref{eq:shnf} satisfying the following 
conditions: (i) the reduced system \eqref{eq:shnfslow} has a singularity at 
$(0,0)$, which is a folded node for $\nu > 0$, (ii) exactly one equilibrium 
point $p_{eq}$ exists in the full system 
with $\nu = O(\varepsilon)$, with a pair of complex conjugate eigenvalues $\rho 
\pm i \omega $ and one real eigenvalue $\lambda$, (iii) the stable manifold 
$W^s$ of $p_{eq}$ is one-dimensional and the unstable manifold $W^u$ of 
$p_{eq}$ 
is two-dimensional ($\lambda < 0 < \rho$). Our notation for $W^s$ and 
$W^u$ hides their dependence on the parameter values. 

We comment on requirements (i) and (ii) above. When 
studying the singular limit $\eps \to 0$ of system \eqref{eq:shnf}, 
condition (ii) requires that $\nu$ must also change, so that the limiting 
system has a folded saddle-node, type II. 

In this regime, small-amplitude oscillations may be due to intersections 
of the attracting and repelling slow manifolds as they twist around each other 
near a folded singularity 
\cite{desroches2012siamreview, wechselberger2005siads} 
or to the spiraling of trajectories near the unstable manifold of the 
equilibrium or both. 
We find that the homoclinic orbits we seek pass through 
the twisting region, so that the interactions of $W^s$ and $W^u$ with the 
slow manifolds of the system play a significant role in their existence. In 
particular, the homoclinic orbits we locate contain segments that lie close to 
the intersection of $W^u$ with the repelling slow manifold $S^r_{\varepsilon}$, 
similar to the homoclinic orbits that form the traveling wave profiles for the 
FitzHugh-Nagumo equation~\cite{GK2009}, a system with one slow variable and two 
fast variables. The homoclinic orbits also contain segments where $W^s$ lies 
close to $S^{a+}_{\varepsilon}$.

To focus upon trajectories that pass through a 
twist region before encountering the equilibrium, we use the position of 
the equilibrium as an alternative to the parameter $\nu$. This position is 
given 
by $p_{eq} = (x_{eq}, x_{eq}^2 + x_{xeq}^3,x_{eq})$, so $\nu = -x_{eq}[a + b 
x_{eq} (x_{eq}+1) + c]$ and the family can be parameterized by $(x_{eq},a,b,c)$ 
instead of $(\nu,a,b,c)$. With this new parameterization, the equilibrium point 
remains fixed as the parameters $a$, $b$, and $c$ are varied, while the 
twist region is a small neighborhood of the origin. In the rest of this 
paper, we set $x_{eq} = -0.03$ whenever we study system \eqref{eq:shnf} 
with $\eps > 0$, except in the concluding remarks where we let the distance 
between $p_{eq}$ and $L_0$ vary with $\eps$.
 
\section{The shooting procedure} \label{sec:shooting}

\begin{figure}[!ht]
\centering
\includegraphics[width=4.5in]{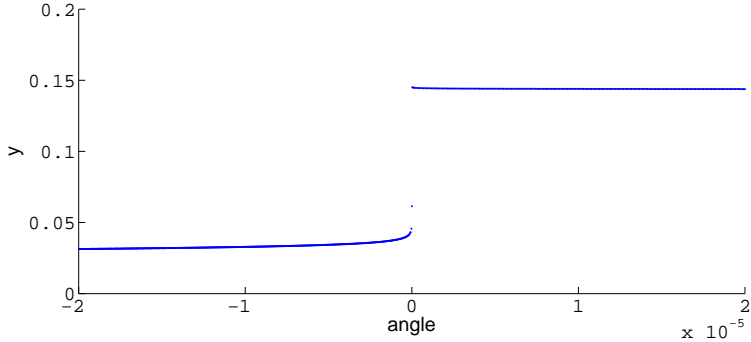}
\caption{Height ($y$ value) of forward trajectories in $W^u$ crossing the 
surface $\{x = -0.05\}$ as a function of angle. An angle of 0 corresponds to a 
trajectory of height $y = 0.1$ lying in $S^r_{\varepsilon} \cap W^u$.
The parameter set is 
$(\varepsilon,a,b,c) \approx (0.01,-4.416165, 2.891404, 5.725663)$. }
\label{fig:heightvsangle}
\end{figure}

The boundary value algorithm HOMCONT~\cite{homcont} was created within the
package AUTO \cite{doedel1996} to compute homoclinic orbits with a collocation 
procedure. Nonetheless, we have been unsuccessful in using AUTO or MATCONT 
\cite{dhooge2003} to locate a Shilnikov homoclinic orbit in \eqref{eq:shnf} 
when 
$\varepsilon \ll 1$. The stiffness of the vector field appears to prevent 
convergence to a homoclinic solution even when very large numbers of 
collocation 
points are used. Shooting algorithms are also problematic since trajectories of 
$W^u$ diverge rapidly from each other near a canard segment of the homoclinic 
orbit. Thus, parameterizing trajectories in $W^u$ by an angular variable and 
varying the angle of an initial point is not well-suited to locating a 
Shilnikov homoclinic orbit because trajectories in $W^u$ are extremely 
sensitive 
to the angle of an initial condition near $p_{eq}$ (see Fig. 
\ref{fig:heightvsangle}). This issue suggests a different shooting strategy 
than 
varying this angle. Instead, we define an angular variable, $\theta$ that 
parametrizes trajectories in $W^u$ smoothly (including as a function of the 
parameters) and regard it as an additional parameter for the system. So $W^u = 
\cup_{\theta} W^u_{\theta}$. The shooting procedure can then fix $\theta$ and 
use another parameter in the search for a homoclinic orbit that contains 
$W^u_{\theta}$.

Our extended family has the six dimensional parameter space $\bar{P}$ with 
coordinates $(\alpha,\theta)$. The homoclinic submanifold of the extended 
parameter space persists as a codimension-one object~\cite{guckenheimer1983}. 
To 
obtain defining equations $\psi:\bar{P} \to \Sigma$ for this manifold, we 
choose 
the surface of section $\Sigma$ defined by $z=0$, set $s^{\alpha}$ to be the 
{\it first} intersection (in backward time) of $W^s$ with $\Sigma$, 
$u^{\alpha}_{\theta}$ to be the {\it first} intersection (in forward time) 
with decreasing $z$ of $W^u_{\theta}$ with $\Sigma$ and 
$\psi(\alpha,\theta) = s^{\alpha} - u^{\alpha}_{\theta}$. The relation 
$\psi(\alpha,\theta)=0$ defines a four dimensional submanifold $\bar{H}$ of 
$\bar{P}$. The projection of $\bar{H}$ to $P$ is a homoclinic manifold 
$H$ consisting of parameters for which \eqref{eq:shnf} has a homoclinic orbit.

Approximations to ${s}^{\alpha}$ and ${u}^{\alpha}_\theta$ are obtained by 
numerically integrating trajectories with initial conditions in the linear 
stable and unstable subspaces of $p_{eq}$. Denoting these approximations by 
$\tilde{s}^{\alpha}$ and $\tilde{u}^{\alpha}_\theta$, the formula 
$\tilde{\psi}(\alpha,\theta) = \tilde{s}^{\alpha} - \tilde{u}^{\alpha}_0$ 
approximates the defining equations. Previous studies \cite{beyn1990contbif, 
beyn1990imajnumanal, lin1990prse,schecter1995imajnumanal} analyze the 
convergence of the solutions $\tilde{\alpha}$ of $\tilde{\psi} = 0$ as the 
distance of the initial conditions to $p_{eq}$ tends to $0$. Hyperbolicity of 
the fixed point, which is satisfied by $p_{eq}$, is required for these 
estimates.

\begin{figure}[!ht]
\centering
\includegraphics[width=4.5in]{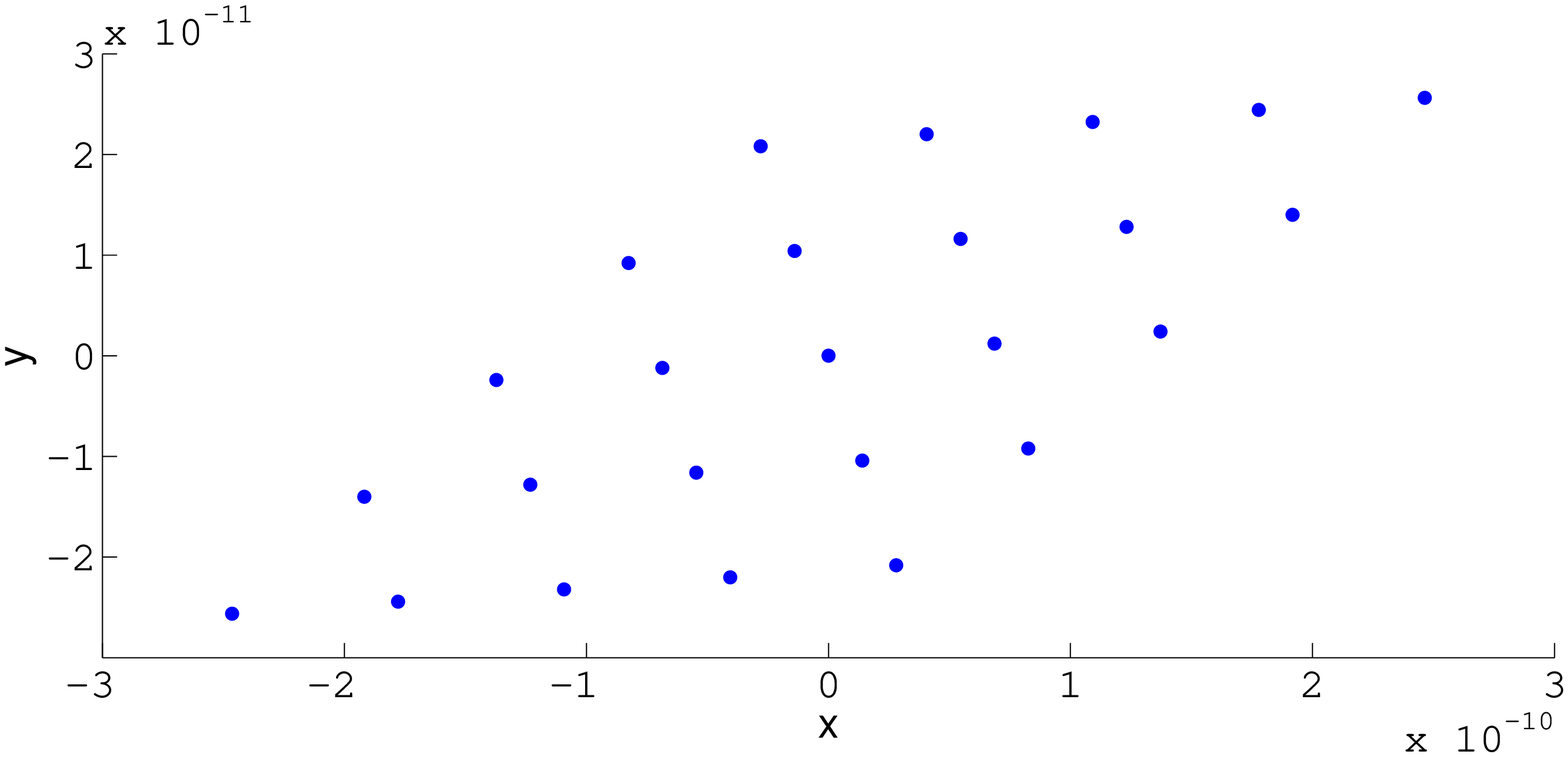}
\caption{The image of a square grid of parameters $(a,b)$ under the 
shooting function $\psi_{\Sigma_0}$ in the surface $\Sigma_0 = \{z = 0\}$. The 
location of the Shilnikov homoclinic orbit is found to be close to the 
parameters $ \tilde{\alpha} \approx (0.01,-0.03,-0.2515348,-1.6508230,1)$. 
Integrations were performed for a square grid of points $(a,b)$ in the 
domain of $\psi_{\Sigma_0}$, specified by 
$a \in [\tilde{a} - 2*10^{-6},\tilde{a} +  2*10^{-6}]$ and $b \in [\tilde{b} - 
2*10^{-6},\tilde{b} +  2*10^{-6}]$.}
\label{fig:homorbit}
\end{figure}

\begin{figure}[!ht]
    \centering
    \subfloat[ ]{{\includegraphics[width=3in, height=1.7in]{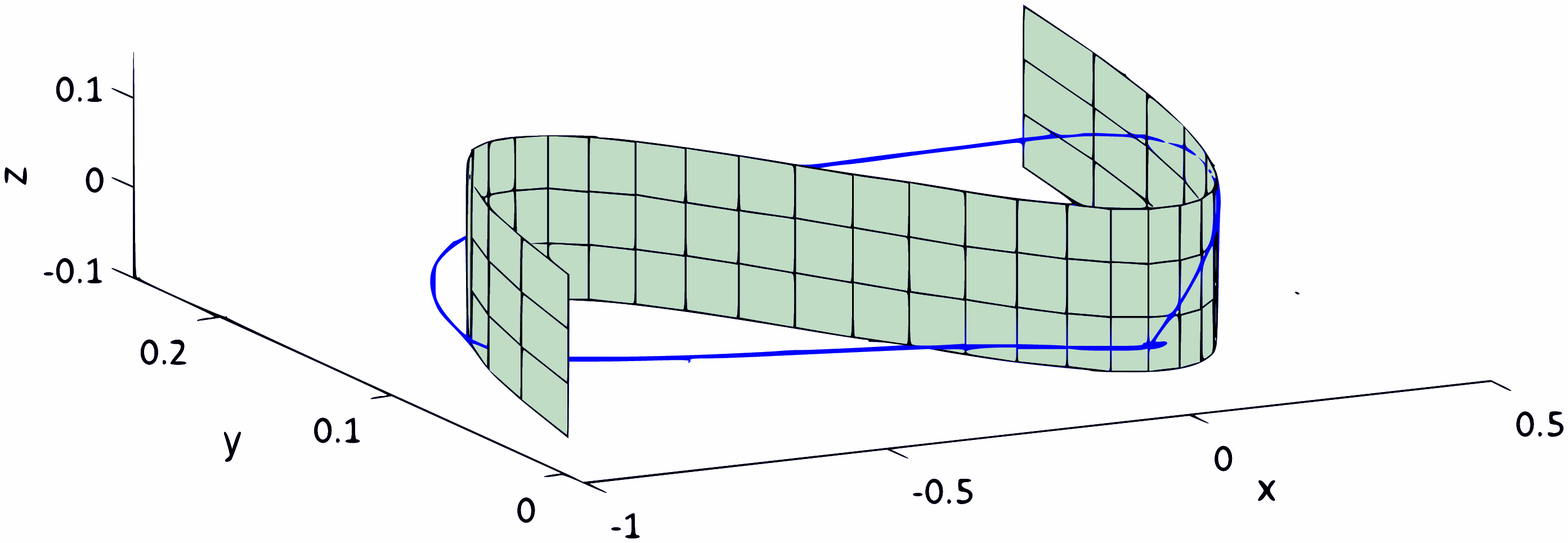} }} 
    \subfloat[ ]{{\includegraphics[width=3in, height=1.7in]{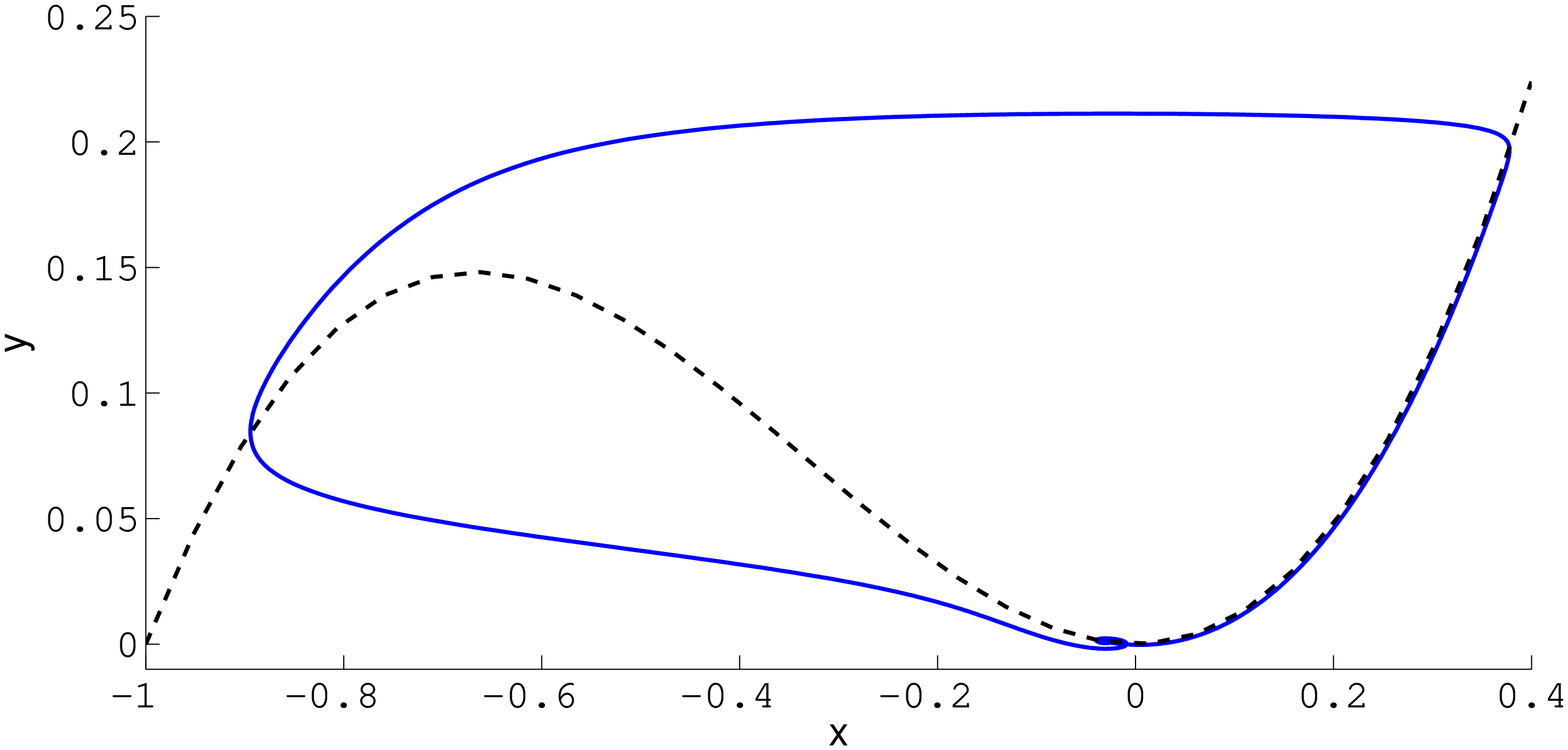} }}
    \caption{Homoclinic orbit (blue curve) to $p_{eq}$ specified by the 
parameters $\theta = 0$ and $\tilde{\alpha}$ in Fig. \ref{fig:homorbit}. The 
critical manifold $C = \{ y = x^2 + x^3\}$ is given by the light green manifold 
in (a) and its $xy$-projection is given by the black dashed curve in (b).}
    \label{fig:homorbitphase}
\end{figure}

We now reduce the number of active parameters by fixing $\varepsilon = 0.01$, 
$\theta = 0$,  $ c = 1$, and $x_{eq} = -0.03$.  Note that our choice of 
$\Sigma$ as the hyperplane $z=0$ is motivated by the complicated dependence of 
$W^s$  on the parameters~\cite{guckenheimer2008chaos}. This section 
slices through a twist region and is close enough to the 
equilibrium point that small changes in $\alpha$ do not produce large jumps in 
$s^{\alpha}$. These choices leave $a$ and $b$ as active parameters to vary in 
$\bar{P}$ to locate an approximate homoclinic orbit by solving the approximate 
defining equations $\tilde{\psi} = 0$.

Fig.~\ref{fig:homorbit} illustrates the regularity of the defining equations 
$\tilde{\psi} = 0$ on a small rectangle $A\subset R^2$ of points in the 
space of active parameters $(a,b)$. Blue dots represent the images of a $5 
\times 5$ lattice of points in $A$ under the shooting function $\tilde{\psi}$. 
The data indicate that $\tilde{\psi}$ is close to affine and regular on $A$ and 
that its image contains the point $(0,0)$, implying that there exists 
$(\tilde{a},\tilde{b}) \in A$ with $\tilde{\psi}(\tilde{\alpha}) = 0$, where 
$\tilde{\alpha}$ is the parameter set with second and third components given by 
$(\tilde{a},\tilde{b})$. This in turn implies the existence of a Shilnikov 
homoclinic orbit, depicted by the blue curve in figure \ref{fig:homorbitphase}.

\section{Continuation of the homoclinic orbit}  \label{sec:continuation}
\begin{figure}[!ht]
\centering
\includegraphics[width=4.5in]{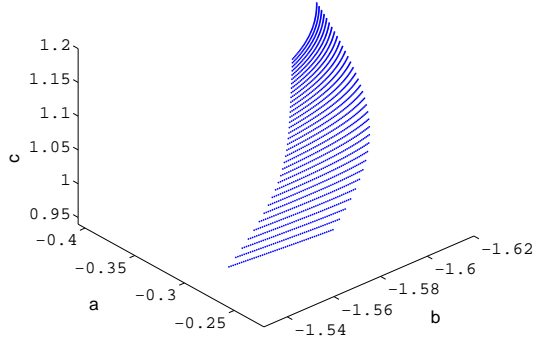}
\caption{A portion of $H$ created with continuation. The portion was projected 
from $(\theta, a, b, c)$ space to $(a,b,c)$ space. The apparent curves show the 
results of continuation in the direction of one of the two nullvectors of $J = 
D\xi$. }
\label{fig:hompatch}
\end{figure}

\begin{figure}[!ht]
    \centering
   \subfloat[ ]{{\includegraphics[width=3.2in, height=2in]{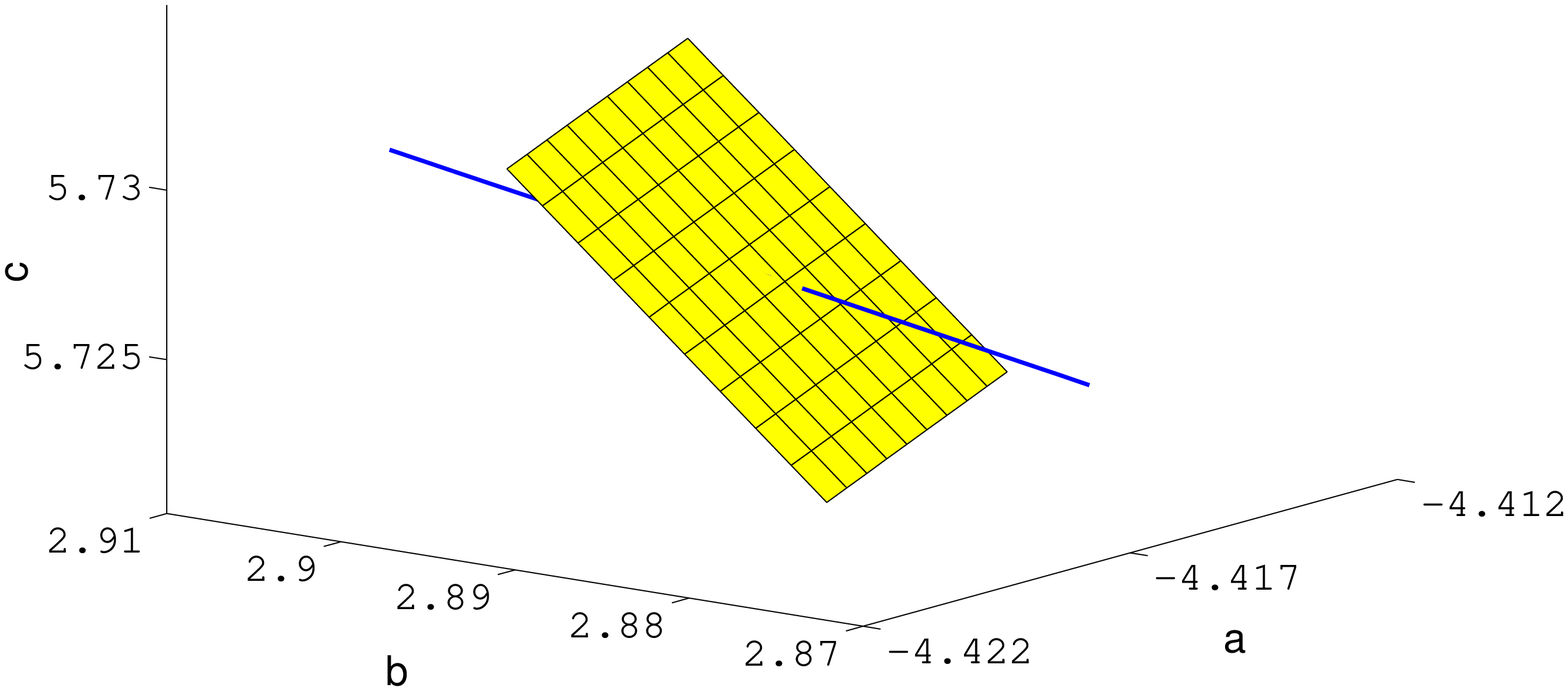} }}
    \subfloat[ ]{{\includegraphics[width=2.5in, height=2in, trim = 5mm 0mm 0mm 
0mm,clip]{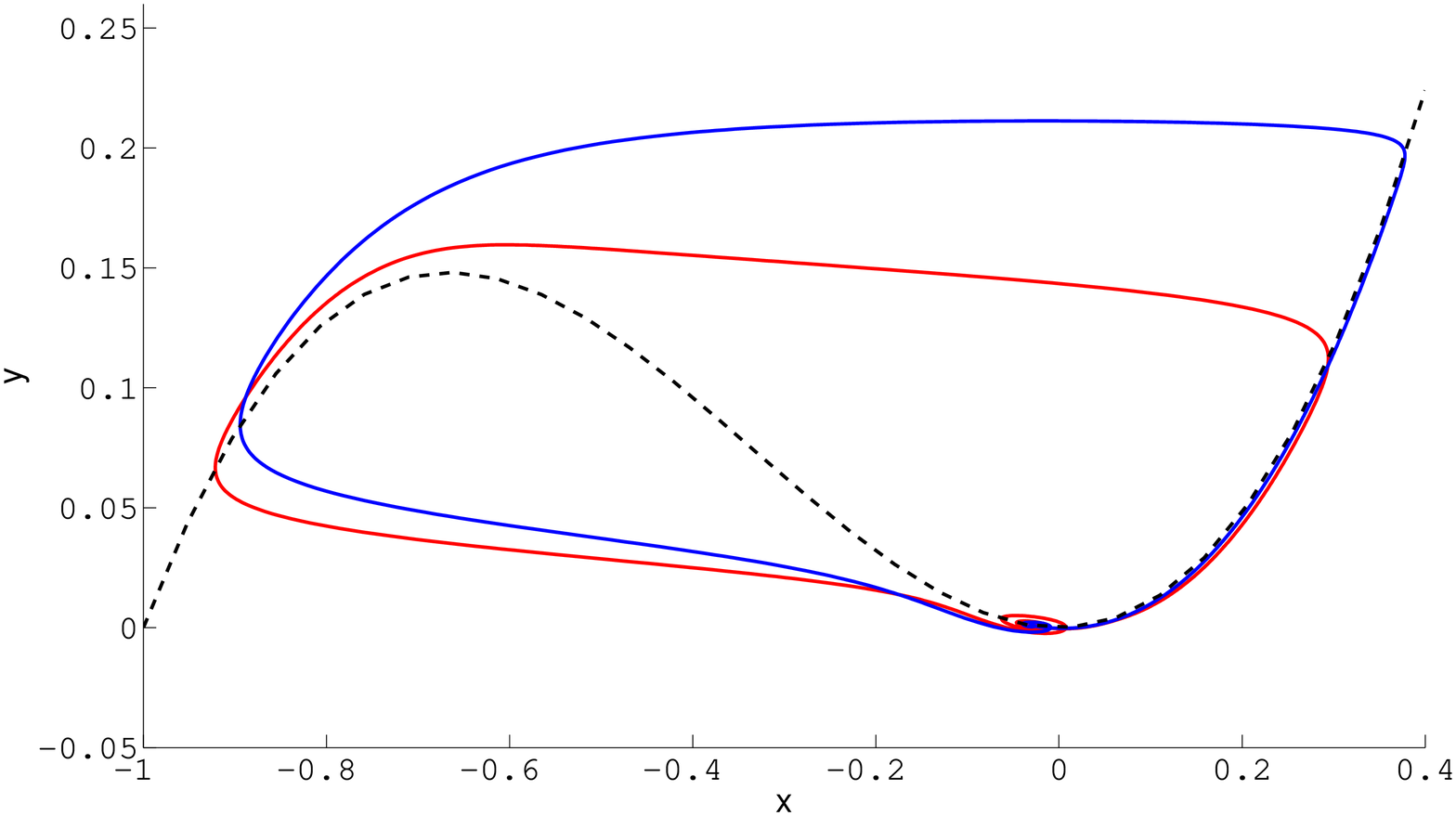} }}
    \caption{(a) Continuation of a curve (blue) on $H_0'$ through a local patch 
of the Koper manifold (yellow)  defined by Eq. \ref{kopeqn}. (b) A 
comparison of the homoclinic orbit defined by the parameter set 
$\tilde{\alpha}$ (blue curve) and the homoclinic orbit lying in the Koper 
subfamily  (red curve) obtained via the continuation in (a), with parameter set 
$\beta = (\varepsilon_{\beta},a_{\beta},b_{\beta},c_{\beta}) \approx 
(0.01,-4.416165, 2.891404, 5.725663)$. The eigenvalues of $p_{eq}$ for the 
Koper homoclinic are $(\rho \pm i \omega, \lambda) \approx (0.790204 \pm  
8.482321 i, -1.576071)$, so $p_{eq}$ satisfies the Shilnikov condition.}
    \label{fig:kophom}
\end{figure}

Continuation algorithms~\cite{kuznetsov1998} are widely used to find curves of 
bifurcations in multi-dimensional parameter spaces. Here, the goal is to find 
an intersection of the homoclinic manifold $H$ with the subfamily of 
\eqref{eq:shnf} that yields the Koper model. The tangent space to $\bar{H}$ is 
the null space of $D\psi$. We estimate $D\psi$ with a central finite difference 
method to provide starting data for continuation calculation of curves on 
$\bar{H}$. These are iterative calculations that use a {\it 
predictor-corrector} algorithm to compute a sequence of parameter values 
$\alpha_j$ on $H$:

\begin{enumerate}
\item {\it Prediction step}: Compute $w^0 = \alpha^j + h v^j$, where $v^j \in 
{\rm ~null}(D\psi)$ and $h$ is our chosen stepsize.
\item {\it Correction step}: Choose a tolerance $\delta$ and iteratively 
compute\\ $w^{k+1} = w^k - (D\psi)^+(w^k) \psi(w^k)$, where $J^+$ is the {\it 
Moore-Penrose pseudoinverse matrix} of $J$ defined by $J^+ = J^T(JJ^T)^{-1}$. 
\item {\it Stopping criterion}: Stop when $||w^{k+1} - w^{k}|| < \delta$ and 
let $\alpha^{j+1} = w^{k+1}$.   
\end{enumerate}

We fix $\varepsilon  = 0.01$, $x_{eq} = -0.03$ and $\theta$, and then use $c$ 
as the third active parameter in addition to $a$ and $b$ in the continuation 
calculation of a curve on $H$. Fig. \ref{fig:hompatch} shows a computation of a 
patch of $H$ in $(a,b,c)$ space with different curves corresponding to 
different values of $\theta$. Fig. (\ref{fig:kophom}a) shows a transversal 
intersection of one of these curves on $H$ with the Koper manifold given by 
Eq.~\ref{kopeqn}. The affine transformations relating \eqref{eq:shnf} to the 
Koper model will rescale and shift the homoclinic orbit (shown in Fig. 
(\ref{fig:kophom}b)) while preserving its topological structure.

\section{Transversality of invariant manifolds} \label{sec:transversality}

\begin{figure}[!ht]
\centering
\includegraphics[width=6in]{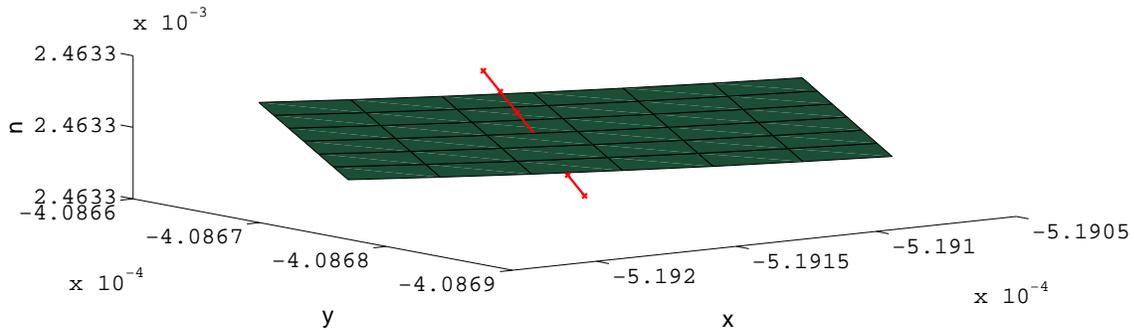}
\caption{Transversal intersection in $\{z = 0\}$ of surface $S$ (green) swept 
out by $S^a_{\varepsilon}$ and curve $C$ (red) swept out by $W^s$ as the 
parameter $a$ is varied. Intersection corresponds to homoclinic orbit defined 
by $\beta$ as given in Fig. \ref{fig:kophom}. Integrations were performed for 
$a \in [a_{\beta} - 3*10^{-5}, a_{\beta} + 3*10^{-5}]$.}
\label{fig:attrslowtostable}
\end{figure}

The homoclinic orbit exists for the full system \eqref{eq:shnf} when a 
trajectory in 
the two dimensional unstable manifold $W^u$ flows along $S^r_\eps$ to a point 
where it jumps to $S^{a-}_\eps$, then flows along $S^{a-}_\eps$ to its fold, 
jumps again to $S^{a+}_\eps$, then flows along this manifold to the folded node 
region where it connects to the one dimensional local stable manifold $W^s$ of 
the equilibrium. This only happens when the parameter values lie in the 
homoclinic 
submanifold $H$ of the parameter space. We can visualize how this happens by 
looking at intersections of $W^s$ and $W^u$ in the cross-section $z=0$ that 
pass through the twist region. On this cross-section, $W^s$ sweeps 
out a curve $C$ and $S^a_{\varepsilon}$ sweeps out a two-dimensional surface 
$S$ 
as the parameter $a$ is varied. Fig. \eqref{fig:attrslowtostable} shows this 
intersection in $(x,y,n)$ space where the local coordinate $n$ is defined via 
$(x,y,a) \cdot  \eta = n$ and $\eta$ is a unit vector normal to a small patch 
of 
$S$. This choice of coordinates increases the angle of intersection that occurs 
in $(x,y,a)$ space. As $a$ varies, the surface swept out by $W^u$ intersects 
the 
curve swept out by $W^s$ 
transversally, demonstrating that the solution of the defining equation for the 
homoclinic orbit is regular.

\begin{figure}[!ht]
\centering

\includegraphics[width=6in]{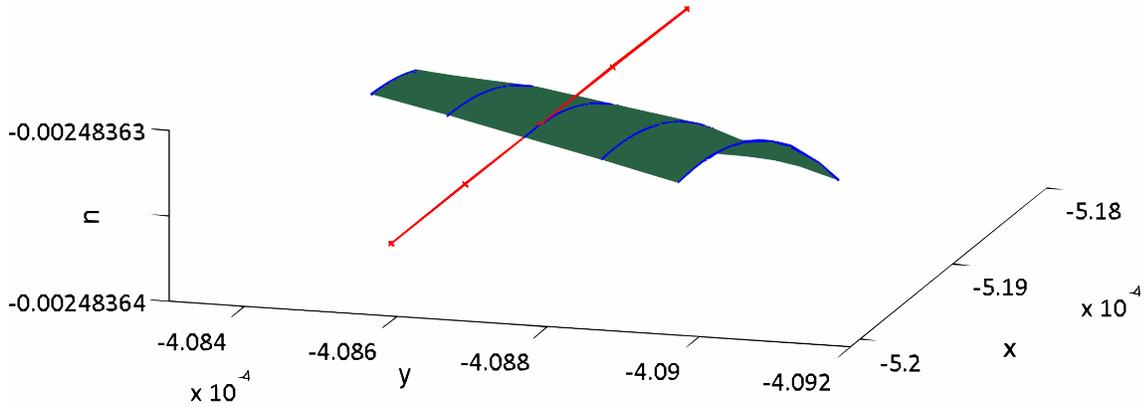}
\caption{Transversal intersection in $\{z = 0\}$ of the surface $S'$ (green) 
swept out by $S^r_{\varepsilon} \cap W^u$ and the curve $C$ (red) swept out by 
$W^s$ as the parameter $a$ is varied. The intersection corresponds to the 
homoclinic orbit defined by $\beta$ as given in Fig. \ref{fig:kophom}. 
Integrations were performed for five equally spaced values of $a \in [a_{\beta} 
- 3*10^{-5}, a_{\beta} + 3*10^{-5}]$. Intersections of $W^u$ (blue curves) for 
different values of $a$ are shown.}
\label{fig:repslowtostable}
\end{figure}

It is difficult to directly compute the relevant portion of $W^u$ since it 
contains canards. 
The two-dimensional subset of $S^r_{\varepsilon}$ we want to compute consists 
of 
trajectories which leave the fixed point $p_{eq}$ along its unstable manifold 
and follow the repelling slow manifold $S^r_{\varepsilon}$ up to some height 
$y$, before finally jumping across to $S^{a-}_{\varepsilon}$. In order to 
locate 
this part of $S^r_{\varepsilon}$, we first integrate backwards a line of 
initial 
conditions at a particular height lying midway between the jump from 
$S^r_{\varepsilon}$ to $S^{a+}_{\varepsilon}$. These backward trajectories jump 
to $S^r_{\varepsilon}$ and flow along it before turning along one of the two 
branches of $W^s$ as they approach the equilibrium point. Since trajectories 
lying in $W^u$ separate those trajectories which follow the two branches of 
$W^s$, we can locate the trajectory $\gamma$ in $S^r_{\varepsilon} \cap W^u$ 
with a bisection method. We then compute an approximation to the strong 
unstable 
manifold of $\gamma$ by integrating forward points on either side of $W^u$ that 
lie close to $\gamma$. This strategy relies on the fact that trajectories 
in $W^u$ approximate leaves of the strong unstable foliation of $S_r$ for the 
canard trajectory in $S^r_{\varepsilon} \cap W^u$. As discussed in Sec. 4, the 
resulting heights of the trajectories are extremely sensitive to the angle as 
illustrated in Fig. (\ref{fig:heightvsangle}). Since the calculation of $W^u$ 
is 
lengthy and indirect, we located the homoclinic orbit parameters by instead 
finding intersections of $W^s$ with $S^{a+}_{\varepsilon}$. Fig. 
(\ref{fig:attrslowtostable}) shows that as we vary the parameter $a$, the 
intersections of $W^s$ with $S^{a+}_\eps$ are similar to those of $W^u$ and 
$W^s$ shown in Fig. (\eqref{fig:repslowtostable}).  As before, this figure 
locates the transversal intersections of $W^u$ and $W^s$ on the surface of 
section specified by $\{z = 0\}$. Trajectories lying on $S^r_{\varepsilon} \cap 
W^u$ sweep out a two-dimensional surface $S'$ and the stable manifold sweeps 
out 
a curve $C$. The objects $S'$ and $C$ intersect transversally  (Fig. 
\ref{fig:repslowtostable}).

\section{Singular homoclinic orbits} \label{sec:singularorbits}

As $\eps \to 0$, trajectories for \eqref{eq:shnf} have singular limits 
consisting of concatenations of fast segments (``jumps'') parallel to the 
$x$-axis and segments that are trajectories of the reduced system on the 
critical manifold. Transitions from slow to fast segments in the singular limit 
trajectories can occur at folds or anywhere along a slow segment on the 
repelling sheet of the critical manifold. The slow trajectory segments are 
contained in invariant manifolds that are approximated by invariant manifolds 
appearing in the singular limit $\eps \to 0$ of the system. These invariant 
manifolds provide a substrate for our theoretical analysis of the homoclinic 
orbits in Section \ref{sec:geomodel}.

\begin{figure}[!ht]
\centering
\includegraphics[width=6in]{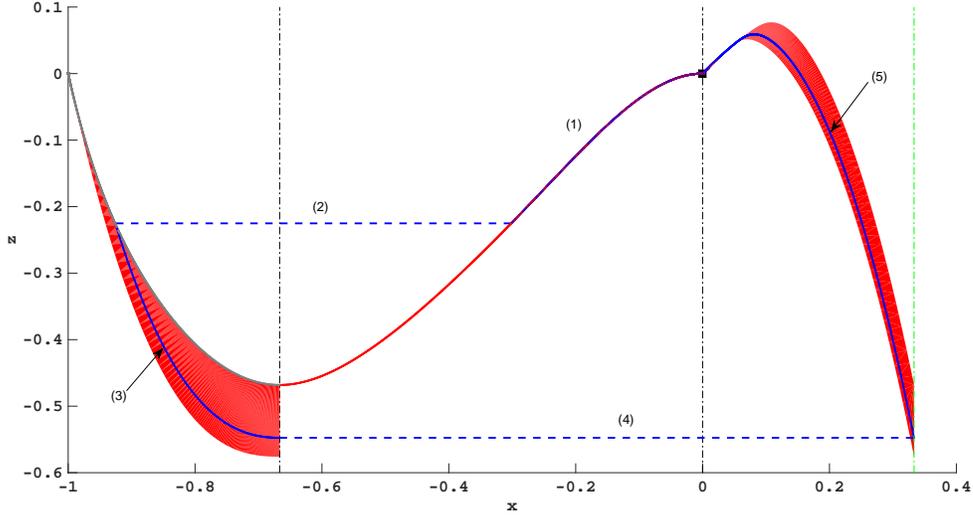}
\caption{ Nonunique singular cycles of the desingularized reduced system 
\eqref{eq:shnfslow}. Recall that stability is reversed on $S^r_0$, the region 
$-2/3 < x < 0$, because of the time reparametrization used to desingularize 
the reduced system. 
A typical cycle in this family contains trajectories on $C$ (blue solid 
curves) and jumps (blue dashed curves): (1) a segment of the unstable 
manifold of the folded saddle-node (0,0), (2) a jump to 
$S^{a-}_0$, (3) a segment of a trajectory flowing forward on $S^{a-}_0$ to the 
left jump curve at $x=-2/3$, (4) a jump to $S^{a+}_0$ ending on the green 
dashed line $\{x = 1/3\}$, (5) a trajectory segment on 
$S^{a+}_0$ connecting the end of this jump to the folded node at (0,0).  
The slow segments of all the singular homoclinic orbits lie in the red regions. 
Parameter set is $(\nu,a,b,c) = (0,a_{\beta},b_{\beta},c_{\beta})$, where 
subscripted parameters are defined  in Fig. \ref{fig:kophom}.}
\label{fig:slowflow}
\end{figure}

When $\eps > 0$, the unstable manifold $W^u$ and the repelling slow manifold 
$S^r_{\eps}$ are 
each two dimensional, so they can intersect transversally along a trajectory. 
Moreover, $W^u$ contains segments aligned 
with the unstable foliation of  $S^r_{\eps}$ beginning where its trajectories 
jump from $S^r_{\eps}$. Tiny variations in initial conditions on $W^u$ 
yield trajectories that turn abruptly at different heights, as illustrated in 
Fig. (\ref{fig:heightvsangle}). This fast portion of $W^u$ turns again to 
follow 
$S^{a-}_{\eps}$ exponentially closely, then jumps to $S^{a+}_{\eps}$ and 
follows this manifold to the twist region. The singular limits of these 
transitions are given by smooth one dimensional maps whose composition with one 
another maps a segment of $W^u$ to a section of $S^{a+}_0$ that passes through 
the folded node. 

We are interested in studying the singular limit of a system for which the 
equilibrium point is a saddle-focus that remains $\eps$-close to the folded 
singularity. These conditions can be satisfied only when the distance from the 
equilibrium point to the fold scales with $\eps$: in the system 
\eqref{eq:shnf}, this implies $\nu = O(\eps)$. In particular, when $\nu = \eps = 
0$ the origin is an equilibrium point of the full system that is a folded 
saddle-node, type II. 
Summarizing, the singular limit of the homoclinic orbits 
can be decomposed as follows: 
\begin{itemize}
 \item 
  An initial segment which lies in $W^u$ within the repelling 
sheet $S^r_0$ of the critical manifold.
 \item
 A jump from $S^r_0$ to the attracting sheet $S^{a-}_0$ of the critical 
manifold.
 \item
 A slow trajectory on $S^{a-}_0$ that ends at the fold at $x=-2/3$.
 \item
 A jump from the fold of $S^{a-}_0$ to $S^{a+}_0$
 \item
 A slow trajectory that follows $S^{a+}_0$ back to the saddle-node 
equilibrium at the origin.
\end{itemize}

We reemphasize here that slow trajectories are defined for the reduced 
system \eqref{eq:shnfslow} at an FSN II bifurcation.

The saddle-node point of this system has a single unstable separatrix and a two 
dimensional stable manifold with boundary. 
This creates a lack of 
uniqueness in the singular limit as shown in Fig. (\ref{fig:slowflow}). In the 
reduced system, trajectories that jump from the unstable manifold to $S^{a-}_0$ 
continue to the fold curve $x = -2/3$ where they jump to points lying in the 
stable manifold of the saddle-node. Jumps from the unstable manifold of the 
saddle-node that occur anywhere on $S^r_0$ yield trajectories that return to 
the saddle-node.
However, when $\eps > 0$, the geometry near the equilibrium 
becomes much more complicated. The equilibrium has a one dimensional 
stable manifold, and this manifold intersects the attracting 
slow manifold $S^{a+}_\eps$ only for parameters lying in a codimension one 
manifold of the parameter space. Moreover, $S^{a+}_\eps$ may twist around
the stable manifold near the location of a folded-node in the reduced system. 
Section \ref{sec:geomodel} investigates the 
persistence of each of the transitions between segments of the singular 
homoclinic orbit as $\eps$ 
becomes positive. Establishing persistence requires transversality hypotheses 
that are formulated there.

\section{Geometry and returns of the twist region} 
\label{sec:returns}

We now turn to a study of the return map to a suitably chosen cross-section 
near the Shilnikov homoclinic orbit we have found.  The classical analysis of 
Shilnikov begins with a homoclinic orbit of a three dimensional vector field 
that has a very special form. First, it is assumed that the vector field is 
linear in a neighborhood $U$ of an equilibrium $p$ and that the eigenvalues 
$\lambda,\rho \pm i \omega$ at $p$ satisfy $|\rho/\lambda| < 1$. Cross-sections 
$\Sigma_1$ and $\Sigma_2$ are chosen in $U$, and the flow map from $\Sigma_1$ 
to $\Sigma_2$ is computed explicitly. The second assumption is that the 
``global return'' from $\Sigma_2$ back to $\Sigma_1$ is an affine map. The 
return map obtained by composing these two flow maps is then proved to have 
hyperbolic invariant sets. Since hyperbolic invariant sets persist under 
perturbation of the vector field, homoclinic orbits of vector fields that do 
not have this special form still have nearby hyperbolic invariant sets. In 
particular, the Shilnikov analysis applies to the homoclinic orbits of 
\eqref{eq:shnf} with the parameter set $\beta$ (as in Fig. \ref{fig:kophom}). 
However, we expect that two aspects of the slow-fast structure of 
\eqref{eq:shnf} may significantly distort the ``standard'' Shilnikov return 
map: (1) the twist region may introduce additional twisting of the 
flow near the homoclinic orbit, and (2) the strong attraction and repulsion to 
the slow manifolds might make the global return map from $\Sigma_2$ to 
$\Sigma_1$ almost singular. We investigate these issues, producing 
modifications of the Shilnikov return map suitable for the homoclinic orbits we 
have located in \eqref{eq:shnf}.

\begin{figure}
\centering
\includegraphics[width=6in]{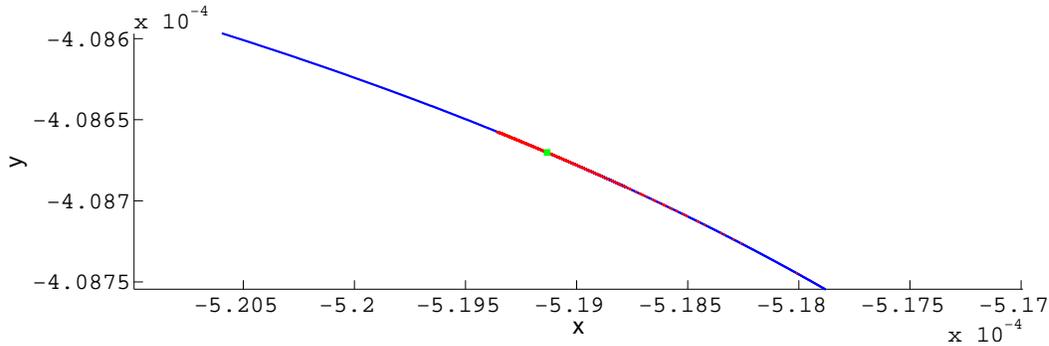}
\caption{A segment of $S^{a+}_{\varepsilon}$ on $\{z = 0\}$ (blue curve) and 
its first return (red points). The first intersection of $W^s$ with $\{z = 0\}$ 
is depicted by the green square. Parameter set is $\beta$.}
\label{fig:spiralreturn}
\end{figure}

\begin{figure}
    \centering
    \subfloat[ ]{{\includegraphics[width=3in, height=2.0in]{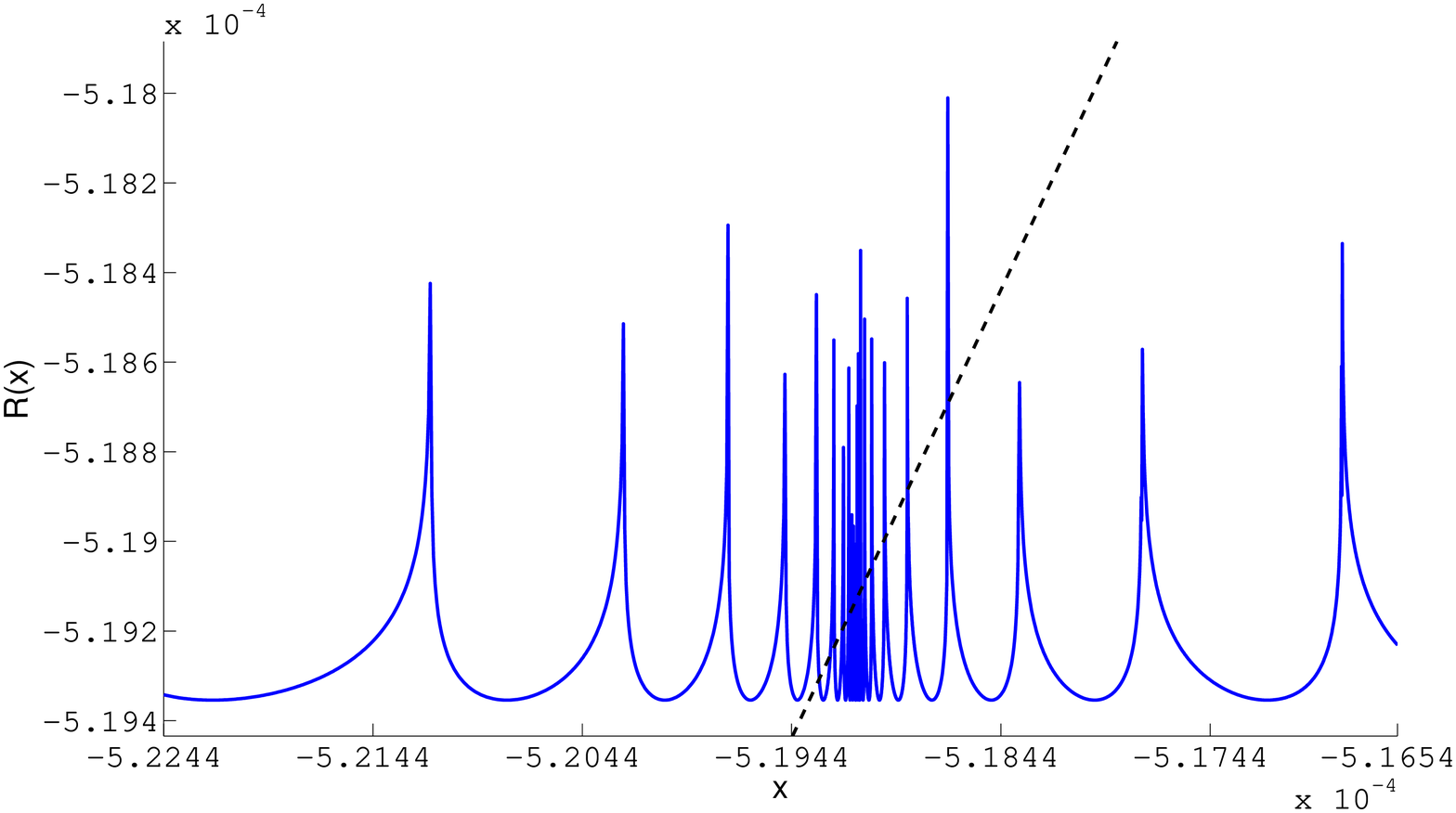} 
}} 
    \subfloat[ ]{{\includegraphics[width=3in, height=2.0in]{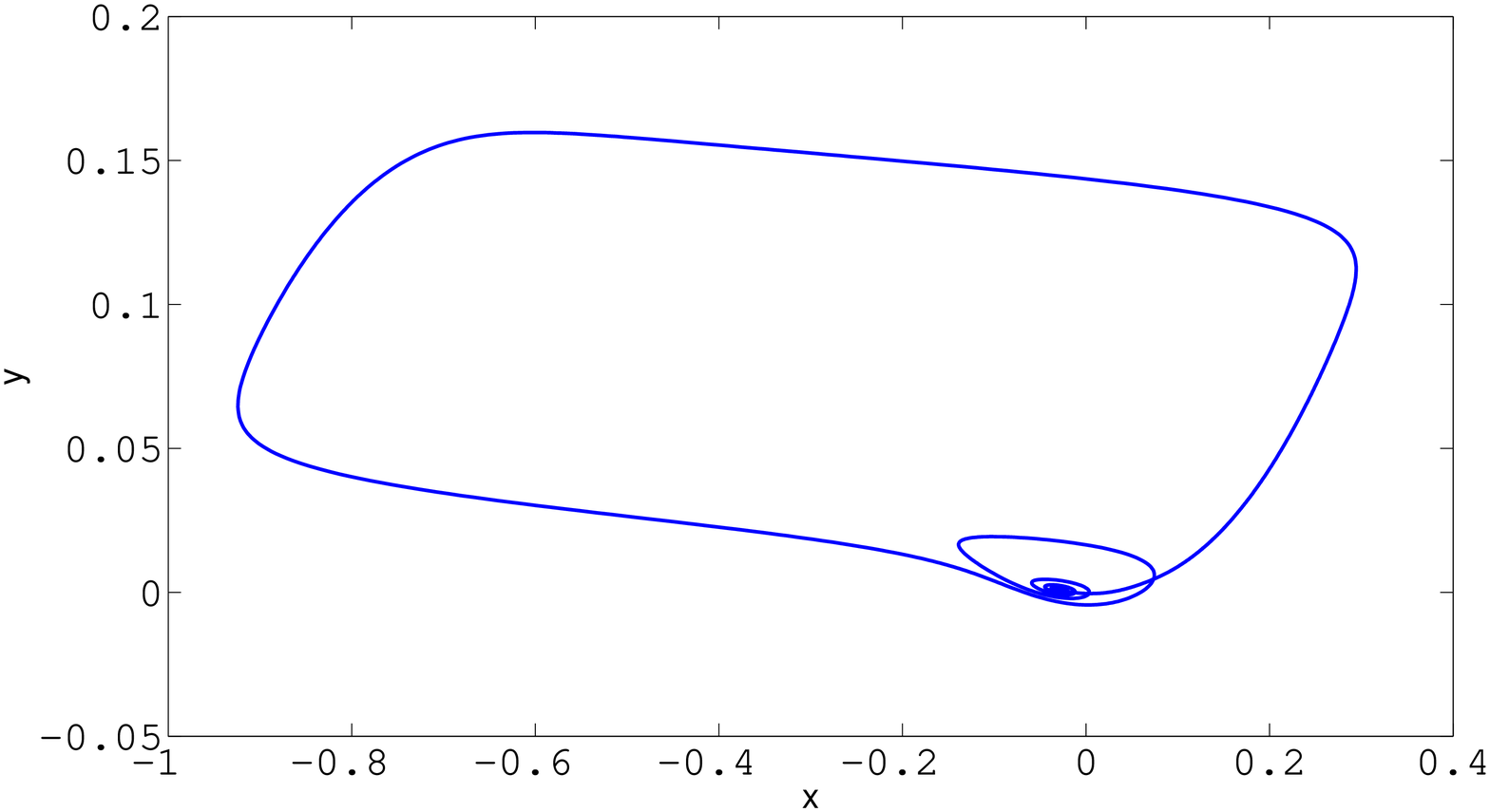} }} 
    \caption{(a) The return map $R$ of points in $I$, where both $(x,y) \in I$ 
and $R(x,y) \in R(I)$ are parametrized by their $x$-coordinates. Points 
$(x,R(x))$ lie on the solid blue curve, and the fixed points that also lie on 
the line $x = R(x)$ (dotted black) belong to periodic orbits that intersect the 
cross section $z=0$ just once. (b) Periodic orbit corresponding to the fixed 
point $p \approx -5.18996*10^{-4}$ of the map $R$.  Parameter set is $\beta$.}
    \label{fig:invreturn}
\end{figure}

We consider the parameter set $\beta$ (as in Fig. \ref{fig:kophom}) and analyze 
the returns of a thin strip $\Sigma$ near $S^{a+}_{\varepsilon} \cap \{z = 
0\}$. 
Exponential contraction of the flow onto $S^{a+}_{\eps}$ suggests that $\Sigma$ 
may be  
mapped into itself in the vicinity of its intersection with the homoclinic 
orbit. Numerical computations suggest that this does happen (Fig. 
(\ref{fig:spiralreturn})). Approximating $\Sigma$ by a small segment $I$ 
parametrized by the $x$-coordinate, the corresponding one-dimensional 
approximation to the return map $R: I \to I$ reveals complicated dynamics (Fig. 
(\ref{fig:invreturn}a)). In particular, we find a sequence of fixed points in 
steep portions of the map $R$ that accumulate at the homoclinic orbit 
intersection. This agrees with previous analyses of homoclinic orbits to 
spiraling equilibrium points that identified a countable number of periodic 
orbits of decreasing Hausdorff distance to the Shilnikov homoclinic orbit 
\cite{GlendinningSparrow,silnikov1965}. The novel behavior here is that these 
fixed points lie in steep portions of the return map whose trajectories contain 
canard segments (Fig. (\ref{fig:invreturn}b)). These periodic orbits may have 
a few
additional twists (small-amplitude oscillations) associated with reinjection 
into the twist region, in addition to the spiraling local to the 
equilibrium point. 

An issue of concern is whether the twisting inside twist regions 
significantly distorts the return map. We explore this issue by comparing the 
geometry we find in the twist region with what is known about the 
twisting of slow manifolds in twist regions. In 
order to make this comparison, we must briefly revisit the slow flow equations. 

In studies of the folded node normal form, Benoit \cite{benoit} 
and Wechselberger \cite{wechselberger2005siads} relate the ratio of 
eigenvalues $\mu = w_1/w_2$ of the folded node in the reduced system to the 
number of intersections of the (extended) attracting and repelling slow 
manifolds. The number $j$ of intersections is estimated by 
the formula $j = 1 + [(\mu-1)/2]$, and this also estimates the maximal number 
of small oscillations of 
trajectories passing through the twist region in the full system. As 
further elucidated by Krupa and Wechselberger~\cite{kw2010}, this estimate 
breaks 
down when the folded node is too close to a folded saddle-node. Here, the 
folded node is too close because $\nu = O(\eps)$ while the results of Krupa and 
Wechselberger~\cite{kw2010} require that $\nu = O(\eps^{1/2})$ or larger. 
Nonetheless, we use the value of $\mu$ as a guide for our numerical 
investigations.  

\begin{figure}
\centering
\includegraphics[width=6in]{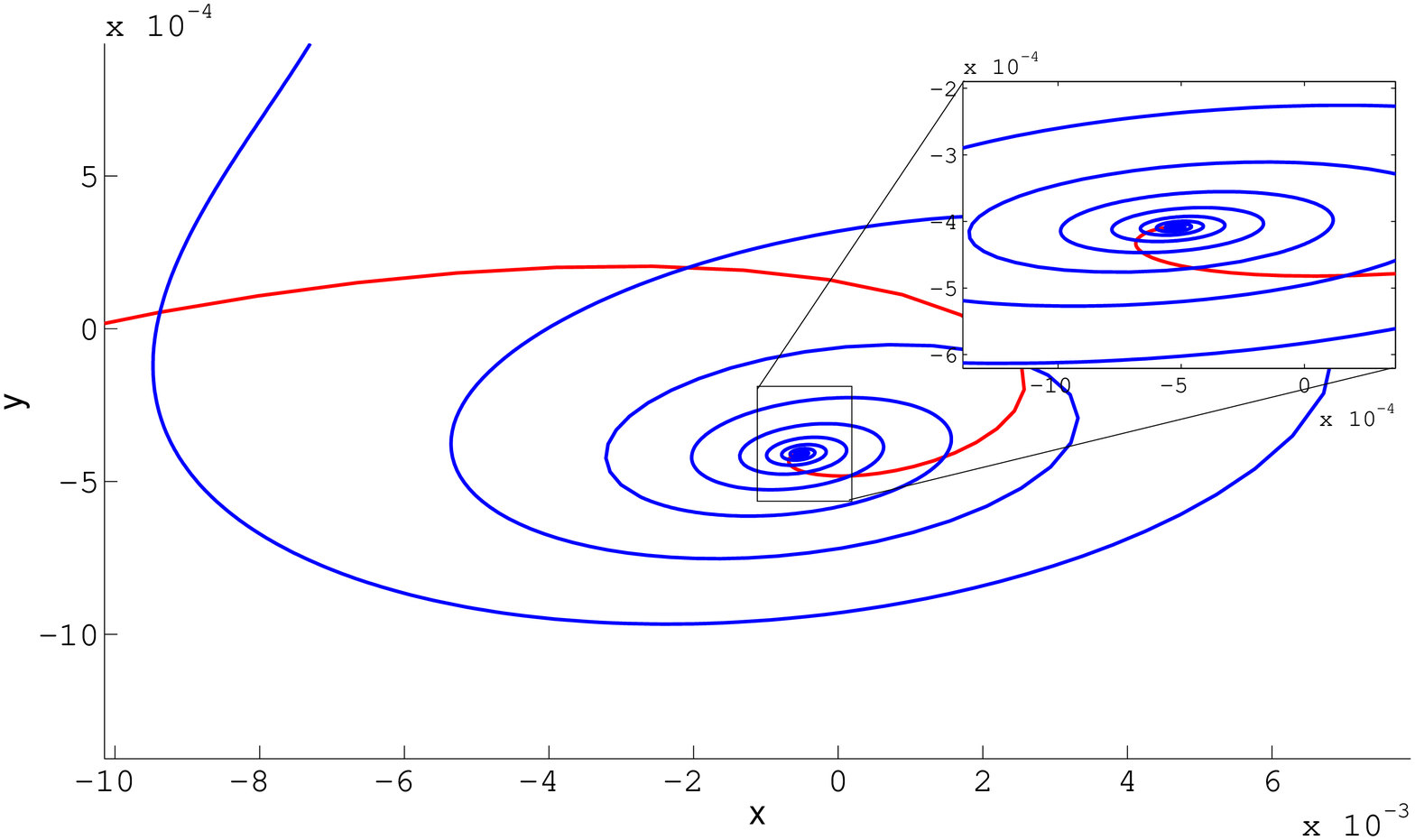}
\caption{Intersection of the extensions of $S^a_{\varepsilon}$ (red) and 
$S^r_{\varepsilon}$ (blue) in the section $\{ z = 0\}$. Inset shows 
magnification of spiraling of $S^r_{\varepsilon}$ near the first intersection 
of $W^s$ with $\{z = 0\}$. Parameter set is $\beta$.}
\label{fig:spirals}
\end{figure}

In the present instance, the folded node of our reduced system equations 
\eqref{eq:shnfslow} has eigenvalues $w_1 \approx -0.920102$ and $w_2 \approx 
-0.0798982$ at the parameter set $\beta$. Thus we estimate $j = 1 + [(11.5159 - 
1)/2] = 6$. However, since the equilibrium point lies in the intersection of 
the extended slow manifolds, $S^r_{\varepsilon}$ has an infinite number of 
turns that yield a countable number of intersections with $S^a_{\varepsilon}$ 
(Fig. (\ref{fig:spirals})). Although the twist region serves to 
strongly contract volumes of the phase space, the equilibrium point produces
the large numbers of small-amplitude oscillations. It seems that the 
twisting inside the twist region does not contribute significantly to the 
geometry of the return map near the homoclinic orbit.

\begin{figure}
\centering
\includegraphics[width=6in]{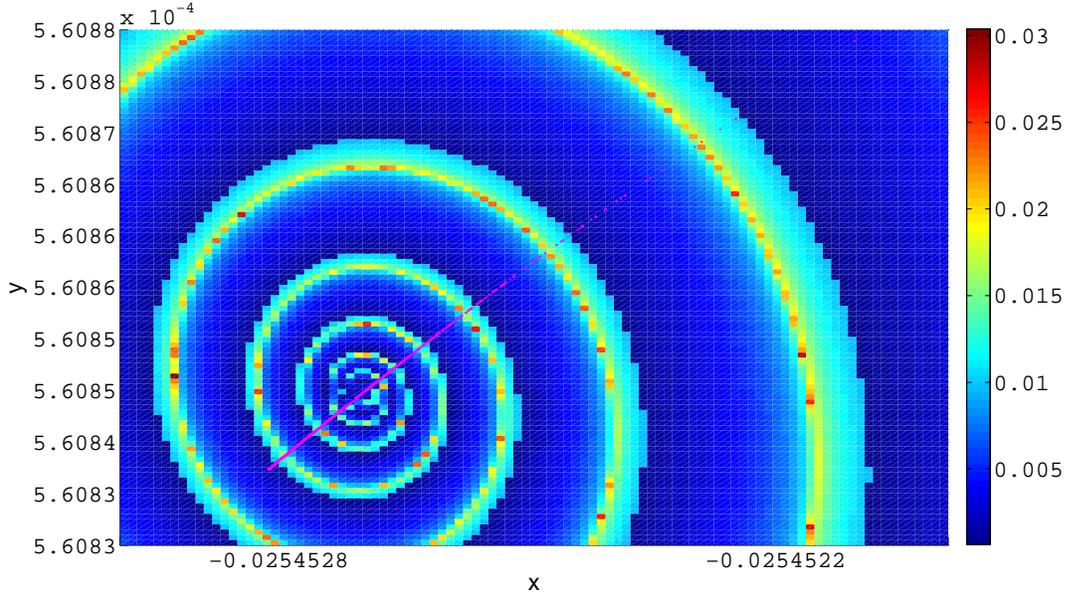}
\caption{Intersection of $I$ (purple) with $S^r_{\varepsilon}$ in the 
cross-section $\{z = -0.025\}$. Color values denote heights ($y$-values) 
attained by trajectories  beginning in a $100 \times 100$ grid of initial 
conditions $(x,y,-0.025)$, where $x \in [-0.025453,-0.025452]$ and $y \in 
[5.6083*10^{-4},5.6088*10^{-4}]$. Stopping conditions are $x = -2/3$ and $y = 
x^2 + x^3 \pm 0.004$. Heights computed from projection of trajectory endpoint 
onto the critical manifold $\{y = x^2 + x^3\}$. Parameter set is $\beta$.}
\label{fig:repspiral}
\end{figure}

Since trajectories beginning in $I$ have canard segments when $I$ intersects 
$S^r_{\varepsilon}$, we examine the resulting distortion by focusing on a 
section closer to the equilibrium point.  Fig. (\ref{fig:repspiral}) shows not 
only that $I$ intersects $S^r_{\varepsilon}$ countably many times near to the 
equilibrium point, but also that canard lengths of forward trajectories are 
organized smoothly in neighborhoods of $S^r_{\varepsilon}$. The property of 
countable intersections is explained by the ``local'' Shilnikov map (in reverse 
time) applied to initial points in $S^r_{\varepsilon}$. Backwards trajectories 
flowing past the equilibrium point spiral very close to $W^s$ by the time the 
trajectory exits a neighborhood of the equilibrium. The distribution of canard 
segment lengths implies that the return map of $I$  stretches and folds subsets 
depending on how the subsets straddle the spiral of the repelling slow manifold.

We get additional insight into the return map using arguments that 
resemble the Exchange Lemma~\cite{joneskopell1994}. This result analyzes the 
Jacobian of a flow map for trajectories that jump from a slow 
manifold of saddle type along its unstable manifold. Consider
the two dimensional system
\begin{eqnarray*}
\varepsilon \dot{x} &=& 1 \\
\dot{y} &=& \lambda y. 
\end{eqnarray*}
\noindent
The flow map of this system from the section $x=0$ to the section $y=1$ is 
given 
by $x(1) =  \frac{-1}{\lambda}\log y(0)$, with derivative 
$\left(\frac{-1}{\lambda y(0)}\right) = \frac{-1}{\lambda} \exp(-\lambda 
x(1))$. 
Note that this derivative grows exponentially with 
the distance a trajectory flows along the slow repelling manifold. For us, 
this is a source of stretching in the global return map of the system 
\eqref{eq:shnf}. We consider trajectories in $W^u$
near $S^r_{\varepsilon}$ that jump to $S^{a-}_{\varepsilon}$. The arrival 
points of these curves can be projected onto  $S^{a-}_{\varepsilon}$ along its 
strong stable foliation. So long as this projection is transverse to the 
trajectories on $S^{a-}_{\varepsilon}$, the stretching that comes from the 
varying jump points on $S^r_{\varepsilon}$ is maintained. Similarly, jumps from 
$S^{a-}_{\varepsilon}$ to $S^{a+}_{\varepsilon}$, projected onto 
$S^{a+}_{\varepsilon}$ along its strong stable foliation, sweep out a curve on 
$S^{a+}_{\varepsilon}$. If this curve is transverse to the 
flow on $S^{a+}_{\varepsilon}$, the stretching in the global return is once 
again maintained.

Of course, there is also fast contraction to the attracting manifolds as well. 
Unless contraction of the flow along $S^{a+}_{\varepsilon}$ in the folded 
node region dominates the stretching originating along 
$S^r_{\varepsilon}$, we can expect that the global return map to be 
approximately a rank one map of large norm for the trajectories that have 
longer 
canards. This is apparent in the spikes of Figure \ref{fig:invreturn}. 

We now verify the claim that stretching is maintained in the global return 
map. Fix two compact, planar cross-sections $\Sigma_0$ and $\Sigma_1$, 
transverse to $W^u$ and $W^s$, respectively. The global return map $R: \Sigma_0 
\to \Sigma_0$ can then be decomposed into the two maps $\varphi: \Sigma_0 \to 
\Sigma_1$ and $\psi: \Sigma_1 \to \Sigma_0$, so that $R = \psi \circ \varphi$. 
Points beginning in $\Sigma_1$ make small oscillations around $W^s$ before 
spiraling out along $W^u$ and hitting $\Sigma_0$. Later, we give an analytical 
approximation of the local map $\psi$. 

We compute $D\varphi(p)$ with central differences, where $p$ is a $10^{-10}$ 
perturbation of the point where the homoclinic orbit intersects $\Sigma_0$. The 
matrix $D\varphi(p)$ has singular values $\sigma_1 \approx 1.53064$ and 
$\sigma_2 \approx 2.66223*10^{-7}$, indicating that the global part of $R$ is 
close to rank one due to strong contraction onto the attracting slow manifolds.

The Jacobian $D\psi$ of the local part of the return map is approximated 
analytically. First, transform coordinates with the real Jordan form $P^{-1} J 
P = J'$, where $J$ is the Jacobian of \eqref{eq:shnf} at $p_{eq}$.  We denote 
transformations of  variables $x$, maps $\zeta$, and subsets $\Sigma$ by primes $x'$, 
$\zeta'$, and $\Sigma'$. Thus, the unstable subspace of $p_{eq}$ becomes
parallel to the 
$x'y'$-plane and the stable subspace becomes parallel to the $z'$-axis. 
Following Shilnikov, $\psi': \Sigma_1' \to \Sigma_0'$ and its derivative 
$D\psi'$ are approximated explicitly from the normal form of the spiraling 
equilibrium point:
\begin{eqnarray}
\psi'(x',y') &=& (r e^{\rho \theta/\omega}, d e^{-\lambda \theta / \omega}), 
\label{eq:locmap}
\end{eqnarray}
where $r = \sqrt{x'^2 + y'^2}$, $\tan \theta = y'/x'$, and $d$ is a small, 
fixed height of $\Sigma_1' $ above $(0,0,0)$.

We recover the Jacobian $D\psi$ by transforming $D \psi'$ to the original 
coordinates on the cross-sections $\Sigma_{0,1}$. By the chain rule, we have 
$DR(p) = D\psi(\varphi(p)) \circ D\varphi (p)$, with eigenvalue magnitudes 
$|\lambda_1| \approx 80166$ and $|\lambda_2| \approx 2*10^{-16}$. We have not 
tried to confirm the relative accuracy of the small eigenvalue, but clearly it 
is very small. Note also that the stretching factor in the local map can be 
shown to become unbounded by picking points approaching the stable manifold on 
$\Sigma_1$ and a sequence of cross-sections $\Sigma_0$ with decreasing heights. 
These points spiral out along the unstable manifold. The number $m$ of turns  
that the trajectory makes before intersecting $\Sigma_0$ determines the 
appropriate solution of the multivalued function $\theta = \arctan(y'/x') + 
2m\pi$ in \eqref{eq:locmap}. The effect on the resulting Jacobian matrix 
$D\psi'$ is multiplication by a diagonal matrix with entries $e^{2m\pi 
\rho/\omega}$ and $e^{-2m\pi\lambda/\omega}$,both of which are larger 
than $1$.

\section{A geometric model of Shilnikov homoclinic orbits} \label{sec:geomodel}

This section abstracts our analysis of the Shilnikov homoclinic orbit in 
the Koper model with a list of geometric conditions 
that are sufficient to prove the existence of such homoclinic orbits in 
slow-fast systems. In the context of this geometric model, the previous 
sections can be regarded, retrospectively, as 
numerical evidence that these conditions are satisfied
along a particular curve of parameter values parametrized by $\eps$ in the 
Koper model. 

The geometric model is formulated in terms of a three dimensional slow-fast 
vector field $X_\eps$ with two slow and one fast variable that depends upon 
additional parameters.

Our first hypothesis is that the reduced 
system $X_0$ (without desingularization) has singular cycles like 
those shown in Fig. \ref{fig:slowflow}:

\textbf{Singular Cycles Hypothesis}:
\begin{itemize}
 \item 
The reduced vector field $X_0$ has an $S$-shaped critical manifold with
sheets $S^{a-}$, $S^r$ and $S^{a+}$ separated by fold curves $L^-$ and $L^+$. 
$S^{a-}$and $S^{a+}$ are attracting while $S^r$ is repelling. The folds are 
generic.
\item 
$X_0$ has a folded saddle-node $p_{eq}$. This point lies on a 
curve of equilibrium points for the full 
system $X_{\varepsilon}$ that are saddle-foci for $\eps > 0 $.
\item
Beginning at $p_{eq}$, the singular cycles consist of 
\begin{enumerate}
 \item
 a segment of the unstable manifold $W^u$ lying entirely in $S^r$,
 \item
 a jump from $W^u$ to $S^{a-}$,
\item
a segment that flows along $S^{a-}$ to $L^-$,
\item
a jump from $L^-$ to $S^{a+}$, 
\item
a segment that flows along $S^{a+}$  back to $p_{eq}$. 
\end{enumerate}
\item
Following jumps from $L^-$ to $S^{a+}$, all of the 
points of $S^{a+}$ from candidates following the first four steps of this 
process form a curve $K$ lying in the basin of attraction of $p_{eq}$. 
\end{itemize}

\textbf{Remark}: In the Koper model and other slow-fast systems with an FSNII 
bifurcation, the equilibrium point is a focal saddle only when it is $O(\eps)$ 
distant from the FSNII point. Thus the second item in this list of hypothesis 
implies that the distance from the equilibrium to the twist region of the 
system scales with $\eps$. System \eqref{eq:shnf} satisfies this scaling 
hypothesis along a curve obtained by setting
$\nu = \varepsilon \bar{\nu}$ and using 
$\bar{\nu}$ as a parameter which is fixed when letting 
$\varepsilon \to 0$.

Proving the persistence of the singular cycles requires additional hypotheses 
that are expressed in terms of transversality. We continue
to refer to the repelling and attracting slow 
manifolds that perturb from the sheets $S^r$ and $S^{a\pm}$ for $\eps > 0$ 
as $S^r_{\eps}$ and $S^{a\pm}_{\eps}$. 

\textbf{Transversality Hypotheses}:
\begin{itemize}
\item
In the singular limit, the image of the jump curve from $W^u$ to $S^{a-}$ is 
transverse to the vector field of the reduced system.
\item
Similarly, the curve $K$ defined above is transverse to the vector field of the 
reduced system on $S^{a+}$. 
 \item 
For $\eps > 0$ small, the unstable manifold of the equilibrium 
point $p_{eq}$ intersects the repelling slow manifold $S^r_{\varepsilon}$ 
transversally in a 
trajectory $\gamma_\eps$.
\item
In the four-dimensional extended phase space that includes a parameter 
$\lambda$, the stable manifolds of $p_{eq}$ sweep out a two-dimensional 
surface as $\lambda$ is varied. This surface intersects the 
three-dimensional attracting slow manifold transversally 
along a trajectory $\beta_{\varepsilon}$. For each $\eps > 0$, 
$\beta_{\varepsilon}$ only exists for a 
particular parameter
value $\lambda = \lambda_h(\varepsilon)$.
\item
The trajectories $\beta_{\varepsilon}$ have a limit as $\eps \to 0$. This 
limit intersects the 
curve $K$ defined above. 
\end{itemize}

\textbf{Remark}: The last item on this list of hypotheses has not been 
investigated thoroughly. Systems with an FSNII bifurcation can be rescaled 
so that the system has a regular limit as $\eps \to 0$
\cite{guckenheimer2008siads,guckenheimer2012siads}. We think that the 
intersections of $W^s$ and $S^{a+}_{\eps}$ that we have analyzed in this paper 
scale nicely with variations of $\eps$ when the remaining parameters are 
suitably scaled, but have little evidence to substantiate this presumption. 
The small amplitude dynamics associated with FSNII bifurcations have not yet 
been studied systematically. 

We now state our main theorem about the geometric model:

\begin{theorem}
 Let $X_\eps$ be a slow-fast vector field with two slow variables and one fast 
variable 
that depends upon an additional parameter $\lambda$. If $X_0$ satisfies 
the 
singular cycle hypothesis, and if $X_\eps$ satisfies the tranversality 
hypotheses, then there is an $\eps_0 > 0$ so that for each $0 < \eps < 
\eps_0$, there is a value of $\lambda = \lambda(\eps)$ for which $X_\eps$ has a 
homoclinic 
orbit.
\end{theorem}

{\it Outline of proof:} Define a cross-section $\Sigma_J$ 
orthogonal to the fast direction in the middle of jumps from $L^-$ to 
$S^{a+}$. Denote by $L^J$ the curve on $\Sigma_J$ that projects onto $L^-$ along 
the fast direction.
Since homoclinic orbits are formed by branches of the stable manifold of 
$p_{eq}$, we prove the theorem by starting at $p_{eq}$ and following its 
stable manifold $W^s$ backward in time. There are intervals of $\lambda$ near 
$\lambda_h(\eps)$ for which the jump points of $W^s$ from 
$S^{a+}_{\eps}$ cross $K$. The fast segments of these jumps intersect 
$\Sigma_J$ in a smooth curve $A$. Projection of $A$ to $S^{a+}_{\eps}$ along 
its fast foliation gives a curve that is close to a trajectory of the reduced 
system on $S^{a+}$. The second transversality hypothesis implies that $L^J$ and 
$A$ are transverse.

Now return to $p_{eq}$ and follow trajectories of its unstable manifold $W^u$ 
until they jump to $S^{a+}_{\eps}$. An exponentially thin wedge of angles in 
$W^u$ gives trajectories that follow $S^r_{\eps}$ for varying distances, 
jumping to  $S^{a-}_{\eps}$ along strong unstable manifolds of $S^r_{\eps}$. 
These trajectories turn to follow $S^{a-}_{\eps}$ where trajectories are 
approximated by trajectories of the reduced system. The first transversality 
hypothesis implies that the width of this strip of trajectories, measured 
orthogonal to the flow direction, will be $O(1)$. When the strip reaches the 
vicinity of $L^-$, it jumps to $S^{a+}_{\eps}$, intersecting $\Sigma_J$ in a 
curve $L^J_{\eps}$. By a classical result of Levinson \cite{Levinson}, $L^J_{\eps}$ is 
$C^1$ close to $L^J$. Consequently, when $\eps > 0$ is small enough, $L^J_{\eps}$ and 
$A$ intersect transversally in $\Sigma_J$ close to a point of $\Sigma_J$ 
lying on a singular cycle. The point $L^J_{\eps} \cap A$ lies on the homoclinic orbit, 
and the theorem is proved.

In addition to proving the existence of the homoclinic orbit, we want to use 
the arguments in Section \ref{sec:returns} to analyze its return map and prove 
that there are chaotic invariant sets nearby. Moreover, these invariant sets 
contain MMOs with unbounded numbers of small amplitude oscillations in their 
signatures:

\begin{theorem}
Let $X_\eps$ be a slow-fast vector field with two slow and one fast variable 
that depends upon an additional parameter $\lambda$. Assume that (1) 
$X_0$ satisfies the singular cycle hypothesis, (2) $X_\eps$ satisfies the 
tranversality hypotheses and (3) that the equilibrium $p_{eq}$ satisfies the 
Shilnikov condition $-\rho/\mu < 1$ for its eigenvalues $\rho \pm i \omega$ and 
$\mu$ when $\eps > 0$. Then, there are chaotic invariant sets in any 
neighborhood of 
the homoclinic orbit of $X_\eps$. These invariant sets include an infinite 
number of periodic orbits that make a single global return around the 
homoclinic orbit. The number of small amplitude oscillations in this set of 
periodic mixed mode oscillations is unbounded. 
\end{theorem}

{\it Outline of proof}: Let $\Sigma_{fn}$ be a cross-section to $X_\eps$ in 
its folded-node region. We 
establish that its return map resembles Figure \ref{fig:invreturn}. Due to the 
strong contraction along $S^{a+}_{\eps}$, this return map will be close to rank 
one 
with image aligned along $I = S^{a+}_{\eps} \cap \Sigma_{fn}$. We study the 
returns of 
a thin strip $\bar{I}$ around $I \subset \Sigma_{fn}$ to $\Sigma_{fn}$. The 
intersection $\xi$ of the stable manifold of $p_{eq}$ with $\Sigma_{fn}$ never 
returns. Shilnikov's original 
analysis \cite{silnikov1965} of the local flow map establishes (1) that points 
that approach $\xi$ make increasing numbers of small amplitude 
oscillations along $W^u$ before flowing along $S^r_{\eps}$, and (2) the 
Jacobian of the flow map 
has a direction with strong expansion. The arguments presented already in 
Section \ref{sec:returns} shows that this expanding direction becomes aligned 
with 
the vector field on $S^r_{\eps}$ as trajectories jump from $S^r_{\eps}$ to 
$S^{a-}_{\eps}$. Arriving 
at $S^{a-}_{\eps}$, our transversality hypotheses imply that the expanding 
direction 
retains a component transverse to the slow flow on $S^{a-}_{\eps}$. The 
transversality 
hypotheses also imply that expansion transverse to the slow flow on 
$S^{a+}_{\eps}$ is 
preserved following the jump from $S^{a-}_{\eps}$ to $S^{a+}_{\eps}$. When 
flowing along 
$S^{a+}_{\eps}$, strong contraction compresses the image of $\bar{I}$ into an 
exponentially thin neighborhood of $S^{a+}_{\eps}$. Thus when points of 
$\bar{I}$ 
return to $\Sigma_{fn}$, their Jacobian is nearly rank one but (by (2) 
above) with strong 
expansion along $I$. 

As in Figure \ref{fig:invreturn}, the spiral 
formed as $\bar{I}$ flows past $p_{eq}$ returns to $\Sigma_{fn}$ with monotone 
branches that cut through $\bar{I}$. Each of these branches contains a fixed 
point at the intersection of a periodic MMO with $\Sigma_{fn}$, and the returns 
that remain within a fixed set of $n$ branches constitute hyperbolic invariant 
sets on which the return map is conjugate to the shift map on sequences of $n$ 
symbols.

\section{Concluding Remark} \label{sec:conclusion}

We have given a fairly complete description of a Shilnikov homoclinic orbit in 
the Koper model, and we have formulated abstract hypotheses that imply 
it occurs in a structurally stable bifurcation for sufficiently small 
$\eps_1$. Our numerical 
investigations suggest that these hypotheses are satisfied.

The Koper model is only moderately stiff in the regime we investigated, 
raising 
the question as to whether the qualitative structure of the homoclinic orbits 
remains unchanged as one approaches the singular limit of the system.
In response to this question, we performed a continuation of the singular Hopf normal form
homoclinic orbit in Fig. (\ref{fig:kophom})  along a parametric path satisfying 
$\nu= \eps \bar{\nu}$ with $\bar{\nu}$ fixed. On such a path, the distance 
from the 
saddle focus to the FSNII scales with $\eps$. As $\eps$ decreases from $0.01$ to 
approximately $0.003$, the resulting 
picture agrees with our analysis of the reduced system (Eq. 
\eqref{eq:shnfslow}): the Shilnikov orbits become better approximated by 
concatenations of slow trajectories on $C$ with jumps across branches of $C$.

\end{document}